\documentclass[11pt]{amsart}

\usepackage{amssymb}
\usepackage{overpic}
\usepackage{enumitem}   
\usepackage{graphicx} 
\usepackage{amsrefs}
\usepackage{xcolor}
\usepackage[hidelinks]{hyperref}

\graphicspath{{Figures/}} 

\newcommand{\R}{\ensuremath{\mathbb{R}}}

\newtheorem {theorem} {Theorem} 

\newtheorem {remark}{Remark}

\begin{document}
	
\title[Piecewise vector fields with hysteresis]{Planar constant piecewise smooth vector fields with large hysteresis}

\author[Tiago Carvalho,  Leonardo Serantola and Bruno de souza Rangel]{Tiago Carvalho$^1$, Leonardo Serantola$^1$ and Bruno de souza Rangel$^2$}

\address{$^1$IBILCE--UNESP, CEP 15054--000, S. J. Rio Preto, S\~ao Paulo, Brazil}
\address{$^2$UFSCar, CEP 13565--905, S\~ao Carlos, S\~ao Paulo, Brazil}
\email{tiago.carvalho1@unesp.br}
\email{l.serantola@unesp.br}
\email{brunodesouzarangel@gmail.com}



\begin{abstract}

Throughout this work, we will carry out a rigorous mathematical analysis of a class of control systems that is widely used in applications but still lacks a consistent theoretical foundation for describing the types of limit sets that may arise from its dynamics. There are applications in which, for example, a treatment for a given disease is administered until the level of diseased cells falls below a prescribed threshold $C_1$. At that point, the treatment is suspended in order to allow the patient’s organism to recover from its side effects. Subsequently, when the level of diseased cells reaches a second threshold $C_2 >C_1$, the treatment is resumed, and the protocol is repeated. To the best of our knowledge, there is not a mathematical classification of such models.

In this paper, we initiate what is intended to become a consistent body of literature aimed at determining the limit sets of such models. We begin with the planar case, in which two linear vector fields are active and two switching boundaries are considered. Naturally, in future developments, control systems in higher dimensions, featuring additional vector fields and more general switching manifolds, should also be considered.
    
\end{abstract}

\maketitle

\section{Introduction}

In this work, we undertake what is, to the best of our knowledge, the first mathematically formal and rigorous study in order to classify the limit sets of a class of control systems that is widely used in applications. Indeed, it is not difficult to envisage problems arising in medicine, agriculture, finance, engineering, chemistry, biology, physics, among other fields, in which such systems naturally emerge.

Here, we will consider two stages (each one modeled by one vector field) and a protocol in the following sense: given an initial condition, the vector field $1$  is considered until the level of a chosen variable $y$  falls below a prescribed threshold $C_1$. At that point, it occurs a change and we consider the vector field $2$ instead of vector field $1$. Subsequently, when $y$ reaches a second threshold $C_2 >C_1$, it occurs a new change and we consider the vector field $1$ again  and the protocol is repeated. To the best of our knowledge, there is not a mathematical classification of such models. Medical applications of this kind of problem are obtained in \cites{Nature-Gatenby-2017, eLife-Gatenby-2022, CarCunEuzFlo-RWA-2025} among others.

When $C_1=C_2$ there is a extensive literature concerning this kind of problems (Filippov systems, for example), but with the deliberate/intentional \textbf{large hysteresis} proposed here, the existing literature is scarce or, in some cases, virtually nonexistent. The fact that $C_1 \neq C_2$ permits that, in a medical protocol, for example, the patient go home and a longer interval between the end of one treatment session and the next (should the number of diseased cells begin to increase again). When $C_1=C_2$ this change can occurs in an arbitrarily small interval of time (for example, when there exists a invisible fold-fold singularity $-$ see \cite{CarBuzTei-JMPA-2014}).

In order to introduce the mathematical description concerning the limit sets of this kind of model, we consider that we should start with the  simplest ones that are:  two planar linear vector fields (dimension $2$) with two switching boundaries $L^{\pm} = \{ (x,y) \in \mathbb{R}^2 : \,   y = \pm \mu\}$ and $\mu>0$. Naturally, in future papers, we wil consider control systems in higher dimensions, featuring additional vector fields and more general switching manifolds. Also, the switching manifolds can not be straight lines. 

In our main results, we show that the trajectories may exhibit:
	$(i)$ unbounded dynamics, with at least one coordinate tending to $\pm\infty$;
	$(ii)$ a monotone \textit{zig-zag} behavior confined to the hysteresis band between $L^\pm$, with the $x$-coordinate tending to $\pm\infty$;
	$(iii)$ periodic orbits within the hysteresis band.

The paper is organized as follows. In Section 2, we introduce the  hysteretic switching systems considered throughout the paper, together with the basic definitions and preliminary remarks. In Section 3, we present the proof of the main result and describe the corresponding asymptotic behaviors of the trajectories. In particular, we characterize the existence of periodic orbits, monotone \textit{zig-zag} dynamics within the hysteresis band and unbounded trajectories. Finally, in Section 4, we present the concluding remarks and discuss possible directions for future research.

\section{Statement of the Main Results}

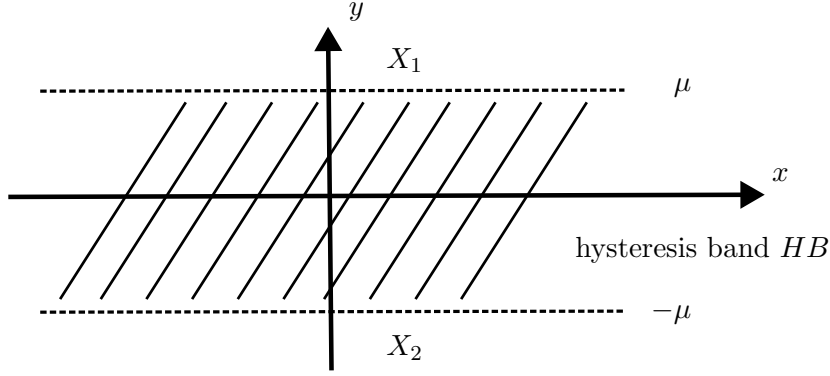
\begin{figure}[ht]
	\begin{center}
		\begin{overpic}[width=10cm]{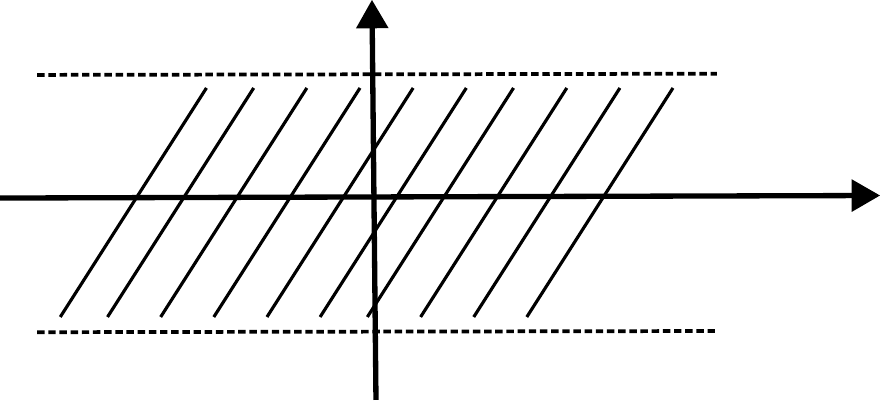}
			\put(88,37){$\mu$}
			\put(85,7){$-\mu$}
			\put(101,25){$x$}
			\put(45,47){$y$}
			\put(50,40){$X_1$}
			\put(50,2){$X_2$}
			\put(75,15){$\mbox{hysteresis band } HB$}
		\end{overpic}
	\end{center}
	\caption{The  configuration studied with an hysteresis band $HB$ between  $y=-\mu$ and $y=\mu$. The vector field $X_1$ acts on the region $\Sigma^+$ and the vector field $X_2$ acts on the region $\Sigma^-$.}\label{figura_0}
\end{figure}

In this section we present the main results concerning the dynamics generated by pairs of vector fields acting  above and below of the large hysteresis band. For $\mu >0$, define the regions
\begin{equation*}
    \Sigma^-= \{ (x,y) \in \mathbb{R}^2 : \,   y \leq - \mu\} \,\, , \,\,  \Sigma^+= \{ (x,y) \in \mathbb{R}^2 : \,   y \geq \mu\} \mbox{ and }
\end{equation*}
\begin{equation}\label{equacao faixa de histerese}
	HB = \{ (x,y) \in \mathbb{R}^2 : \,   - \mu \leq y \leq  \mu\}.
\end{equation}
Throughout the paper we denote by $X_1$ the vector field acting on the region $\Sigma^+$ and by $X_2$ the vector field acting on the region $\Sigma^-$. Both of them are acting in the hysteresis band $HB$ given by \eqref{equacao faixa de histerese}. See Figure~\ref{figura_0}.

\begin{remark}
	\textbf{IMPORTANT CONVENTION:} In this paper we will consider that when an initial condition of a solution belongs to $HB$, then the vector field to be considered is $X_2$. Physically, we are thinking that in $HB$ the ''treatment'' $X_1$ is not applied. An analogous approach can be done considering the opposite case. 
\end{remark}

\begin{remark}
	We denote by
\[
\overline{\mathbb{R}}
:=
\mathbb{R}\cup\{-\infty,+\infty\}
\]
the extended real line. The extended plane is defined as
\[
\overline{\mathbb{R}}^2
:=
\overline{\mathbb{R}}\times\overline{\mathbb{R}}.
\]
    Consider the vector field $X=(X_1,X_2)$ and an arbitrary initial condition $(x_0,y_0)$. Let us define the $\omega$-limit set of the trajectory through $(x_0,y_0)$ by
\[
\omega\!\left(x_0,y_0\right)
=
\left\{
(\omega_1,\omega_2)\in\overline{\mathbb{R}}^2 :
\lim_{n\to\infty}
X\!\left(t_n;x_0,y_0\right)
=
(\omega_1,\omega_2)
\text{ for  }
t_n\to+\infty
\right\}.
\] 
\end{remark}

Throughout the paper, we consider the constant vector fields acting on both sides of the hysteresis band. On the region $\Sigma^+$ we take
\[
	X_1 = (a_1,b_1), \quad a_1,b_1\in\mathbb{R} \quad \textrm{and} \quad a_1^2 + b_1^2 \neq0
\]
and on the region $\Sigma^-$ we take
\[
	X_2 = (a_2,b_2), \quad a_2,b_2\in\mathbb{R} \quad \textrm{and} \quad a_2^2 + b_2^2 \neq0.
\]

	
	
	
		
		
	

\begin{theorem}\label{theorem_1} 
	Consider two constant vector fields
\[  X_1 = ( a_1 ,  b_1  ) \mbox{ and } X_2= (a_2 ,  b_2 ),  \]with $a_1 ,  a_2, b_1,b_2 \in \R $, $a_1^2 + b_1^2 \neq0$, $a_2^2 + b_2^2 \neq0$  and the hysteresis band given by \begin{equation}\label{equacao faixa de histerese}
	HB = \{ (x,y) \in \mathbb{R}^2 : \,   - \mu \leq y \leq  \mu\}, \mbox{ with } \mu>0.
\end{equation}
Given an arbitrary initial condition $(x_0,y_0)$, the  $\omega$-limit set of a trajectory passing through this point is one of the following topological cases:

$\bullet$ $
\{(\pm\infty,\pm\infty)\}$;

$\bullet$ either $
\{(\pm\infty,M)\}$ or $
\{(M,\pm\infty)\}$, with $M$ being a constant. The value of $M$ is determined in the proof of the theorem according to the values of $a_1,a_2,b_1,b_2$;

$\bullet$ periodic orbit restricted to $HB$. Moreover, in this case there exists a continuum of periodic orbit restricted to $HB$;

$\bullet$  either $
\{(\pm\infty,\mbox{ zig-zag })\}$ or $
\{(\mbox{ zig-zag },\pm\infty)\}$, with the zig-zag occurring  inside$HB$.
\end{theorem}

\section{Proof of the Main Results}
	
\subsection{Proof of Theorem \ref{theorem_1}}

\begin{figure}[ht]
	\begin{center}
		\begin{overpic}[width=6cm]{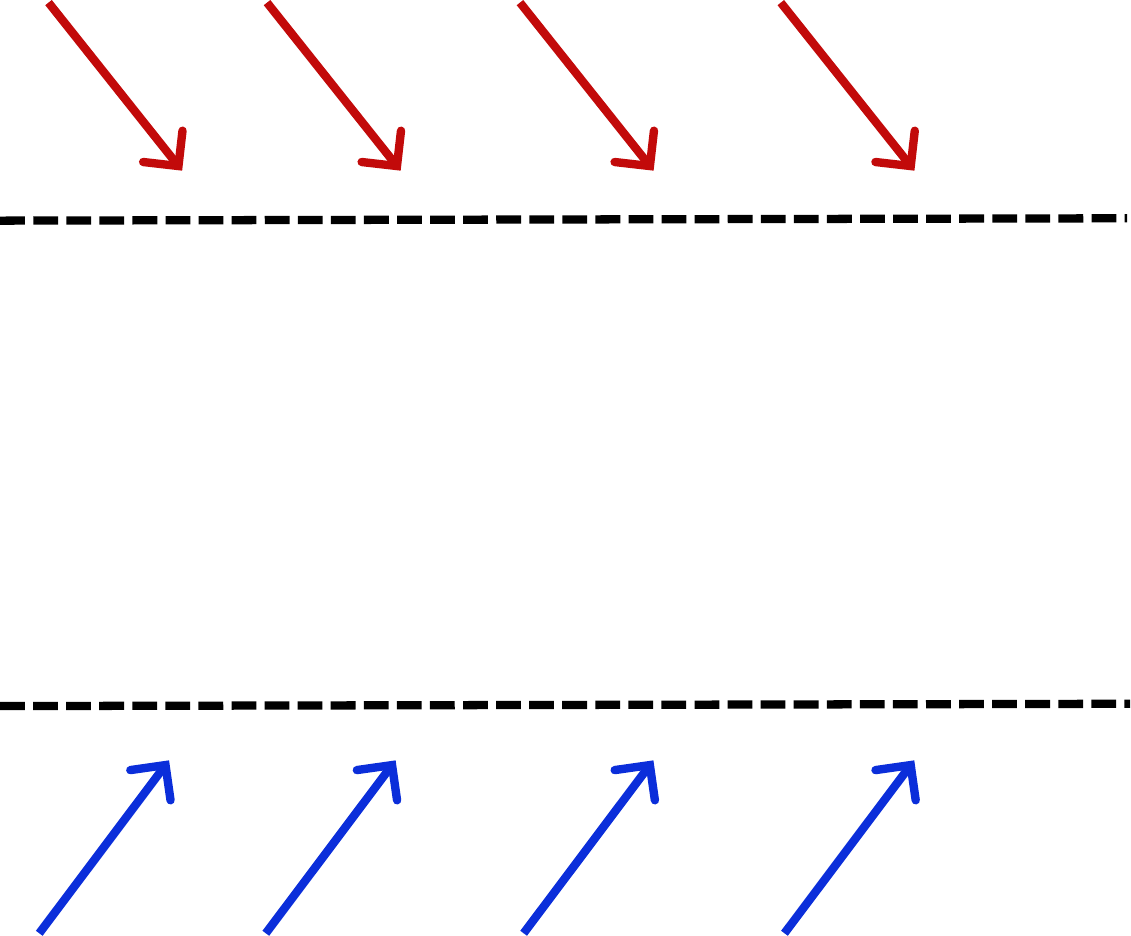}
			\put(106,63){$\mu$}
			\put(103,20){$-\mu$}
			\put(85, 80){$X_1=\left(a_1, b_1\right)$}
			\put(85,3){$X_2=\left(a_2, b_2\right)$}
		\end{overpic}
	\end{center}
\caption{An example of two constant vector fields acting in the upper region and the down region of the hysteresis band with $b_1<0$ and $a_1,\, a_2,\, b_2>0$.}\label{figura_1}
\end{figure}
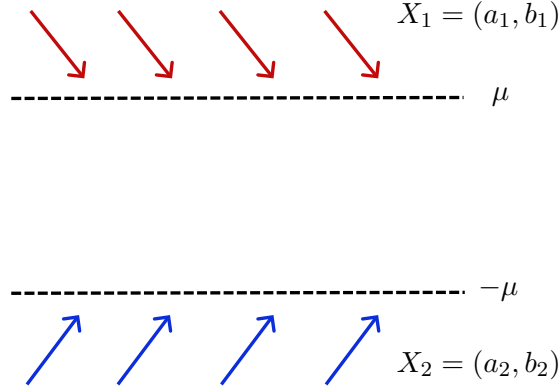

The vector fields in both sides of the hysteresis band are constant.
See Figure \ref{figura_1}.


\vspace{0.5cm}

\begin{figure}[ht]
	\begin{center}
		\begin{overpic}[width=7cm]{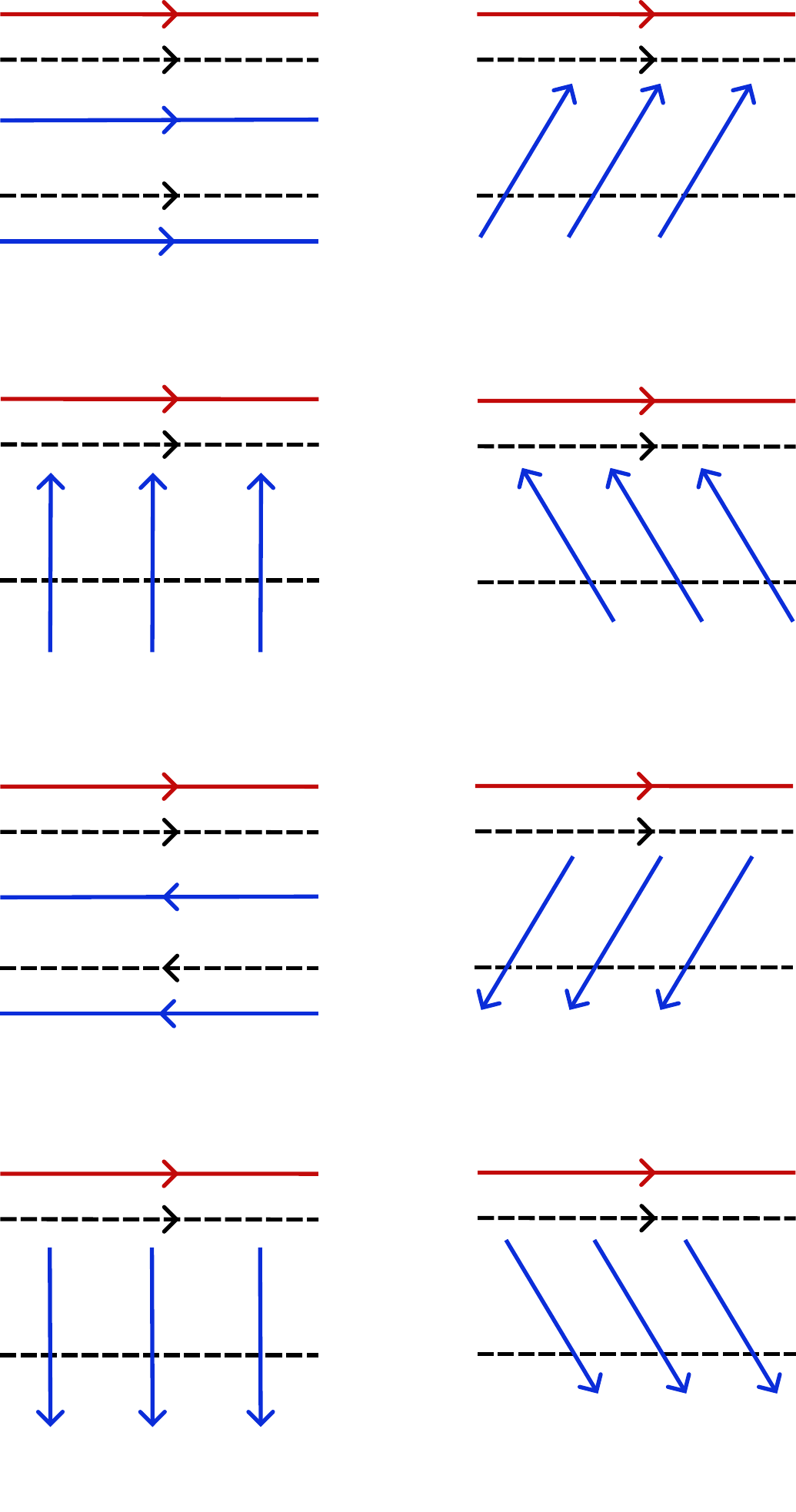}
			\put(6,78){Case $1.1$}
			\put(37,78){Case $1.2$}
			\put(6,52){Case $1.3$}
			\put(37,52){Case $1.4$}
			\put(6,27){Case $1.5$}
			\put(37,27){Case $1.6$}
            \put(6,0){Case $1.7$}
            \put(37,0){Case $1.8$}
		\end{overpic}
	\end{center}
	\caption{Dynamics for conditions $a_1>0$ and $b_1=0$.}\label{figura_1.1}
\end{figure}

\begin{figure}[ht]
	\begin{center}
		\begin{overpic}[width=9cm]{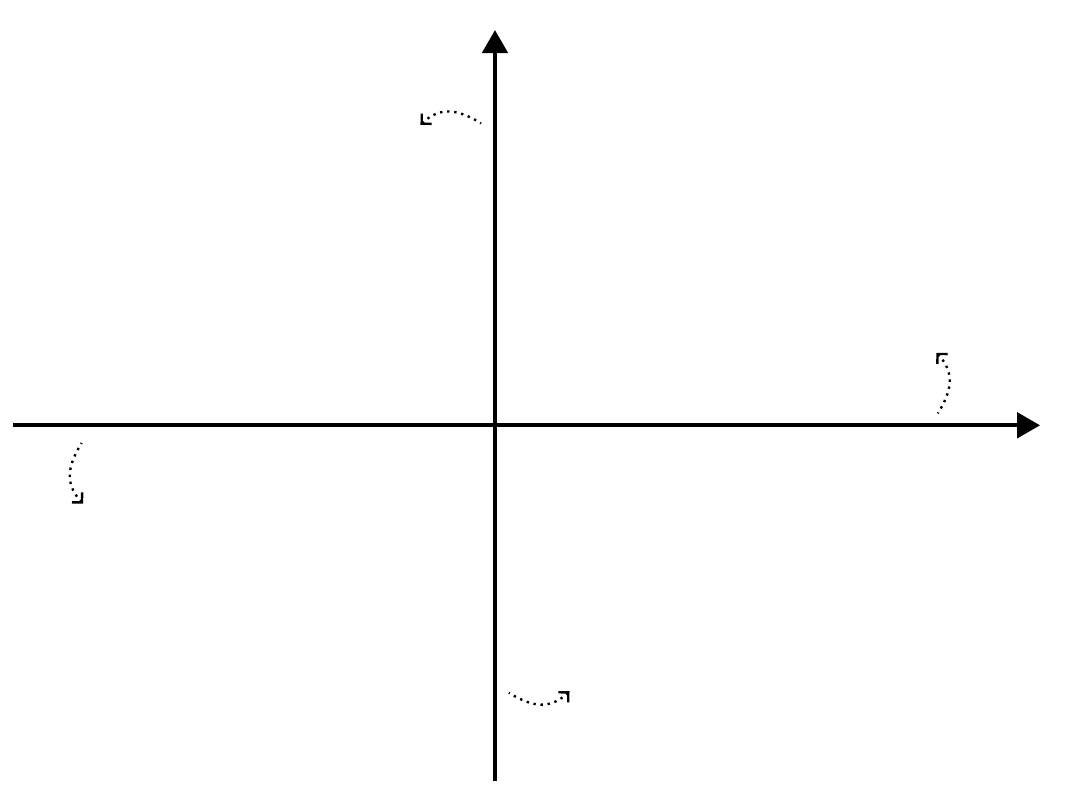}
            \put(97,30){$a_2$}
            \put(48,70){$b_2$}
			\put(75,45){Case $1.1$}
            \put(55,55){Case $1.2$}
            \put(25,55){Case $1.3$}
            \put(5,45){Case $1.4$}
            \put(5,20){Case $1.5$}
            \put(25,10){Case $1.6$}
            \put(55,10){Case $1.7$}
            \put(75,20){Case $1.8$}
		\end{overpic}
	\end{center}
	\caption{Bifurcation diagram in the variables $a_2$ and $b_2$, for $a_1 >0$ and $b_1=0$.}\label{figura_2.1}
\end{figure}

\noindent Here we will analyze the dynamics for all cases where $a_1 >0$ and $b_1=0$. See Figure \ref{figura_1.1} for the phase portraits corresponding to these cases and Figure \ref{figura_2.1} for the bifurcation diagram in the variables $a_2$ and $b_2$.

\vspace{0.5cm}

$\bullet$ \textbf{Case 1.1: $a_1>0, b_1 = 0, a_2 > 0, b_2 = 0$.} Given an arbitrary initial condition $(x_0,y_0)$, the trajectory has $\omega$-limit set $\{(+\infty,y_0)\}$.

$\bullet$ \textbf{Case 1.2: $a_1>0, b_1 = 0, a_2 > 0, b_2 > 0$.} Given an arbitrary initial condition $(x_0,y_0)$, the trajectory has $\omega$-limit set $\{(+\infty,y_0)\}$ when $y_0 > \mu$ and has $\omega$-limit set $\{(+\infty,\mu)\}$ when $y_0 \leq \mu$.

$\bullet$ \textbf{Case 1.3: $a_1>0, b_1 = 0, a_2 = 0, b_2 > 0$.} Given an arbitrary initial condition $(x_0,y_0)$, the trajectory has $\omega$-limit set $\{(+\infty,y_0)\}$ when $y_0 > \mu$ and has $\omega$-limit set $\{(+\infty,\mu)\}$ when $y_0 \leq \mu$.

$\bullet$ \textbf{Case 1.4: $a_1>0, b_1 = 0, a_2 < 0, b_2 > 0$.} Given an arbitrary initial condition $(x_0,y_0)$, the trajectory has $\omega$-limit set $\{(+\infty,y_0)\}$ when $y_0 > \mu$ and has $\omega$-limit set $\{(+\infty,\mu)\}$ when $y_0 \leq \mu$.

$\bullet$ \textbf{Case 1.5: $a_1>0, b_1 = 0, a_2 < 0, b_2 = 0$.} Given an arbitrary initial condition $(x_0,y_0)$, the trajectory has $\omega$-limit set $\{(+\infty,y_0)\}$ when $y_0 > \mu$ and has $\omega$-limit set $\{(-\infty,y_0)\}$ when $y_0 \leq \mu$.

$\bullet$ \textbf{Case 1.6: $a_1>0, b_1 = 0, a_2 < 0, b_2 < 0$.} Given an arbitrary initial condition $(x_0,y_0)$, the trajectory has $\omega$-limit set $\{(+\infty,y_0)\}$ when $y_0 > \mu$ and has $\omega$-limit set $\{(-\infty,-\infty)\}$ when $y_0 \leq \mu$.

$\bullet$ \textbf{Case 1.7: $a_1>0, b_1 = 0, a_2 = 0, b_2 < 0$.} Given an arbitrary initial condition $(x_0,y_0)$, the trajectory has $\omega$-limit set $\{(+\infty,y_0)\}$ when $y_0 > \mu$ and has $\omega$-limit set $\{(x_0,-\infty)\}$ when $y_0 \leq \mu$.

$\bullet$ \textbf{Case 1.8: $a_1>0, b_1 = 0, a_2 > 0, b_2 < 0$.} Given an arbitrary initial condition $(x_0,y_0)$, the trajectory has $\omega$-limit set $\{(+\infty,y_0)\}$ when $y_0 > \mu$ and has $\omega$-limit set $\{(+\infty,-\infty)\}$ when $y_0 \leq \mu$. 

\vspace{0.5cm}

\begin{figure}[ht]
	\begin{center}
		\begin{overpic}[width=7cm]{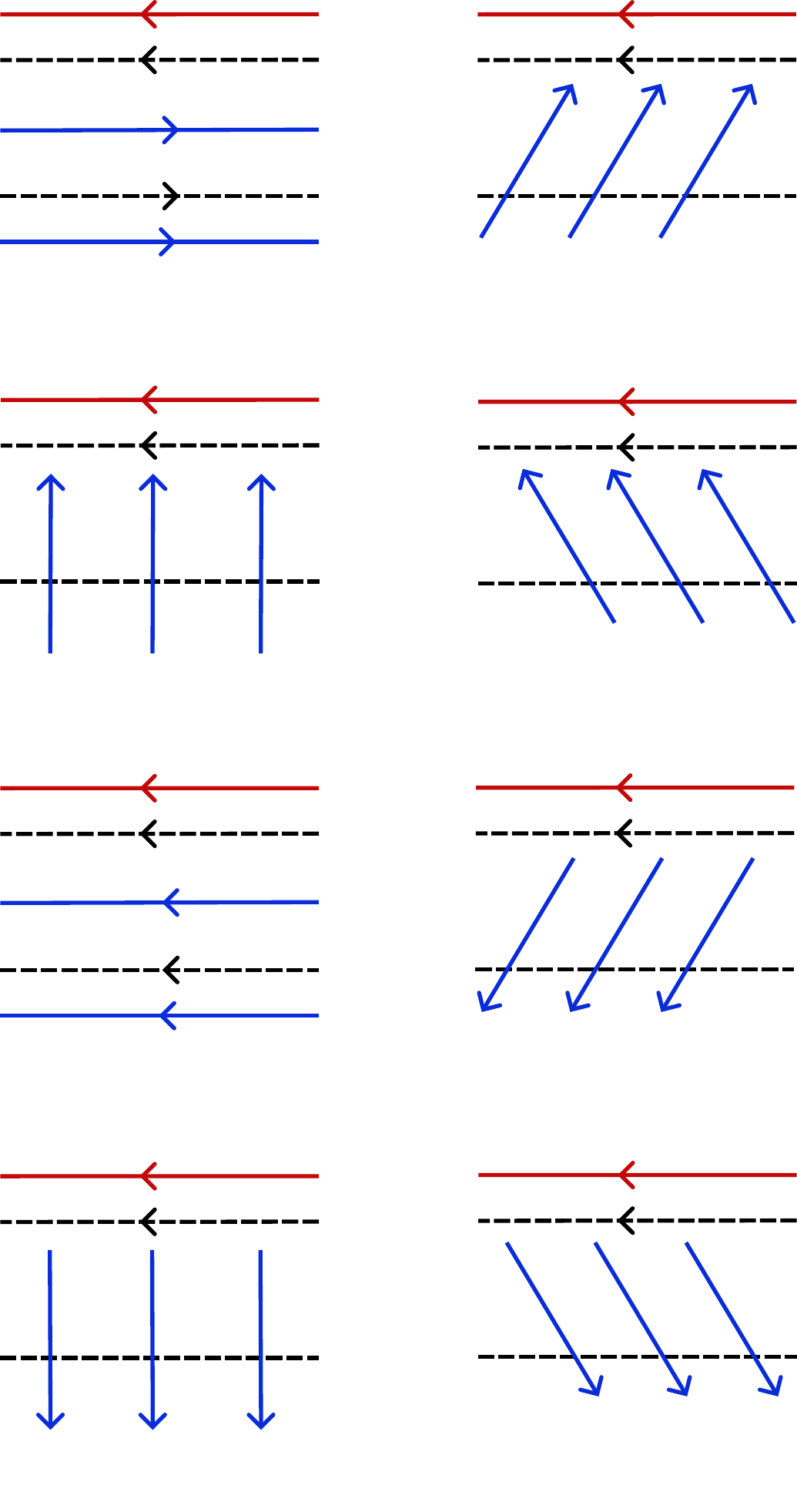}
			\put(6,78){Case $1.9$}
			\put(37,78){Case $1.10$}
			\put(6,52){Case $1.11$}
			\put(37,52){Case $1.12$}
			\put(6,27){Case $1.13$}
			\put(37,27){Case $1.14$}
            \put(6,0){Case $1.15$}
            \put(37,0){Case $1.16$}
		\end{overpic}
	\end{center}
	\caption{Dynamics for conditions $a_1<0$ and $b_1=0$.}\label{figura_1.2}
\end{figure}

\begin{figure}[ht]
	\begin{center}
		\begin{overpic}[width=9cm]{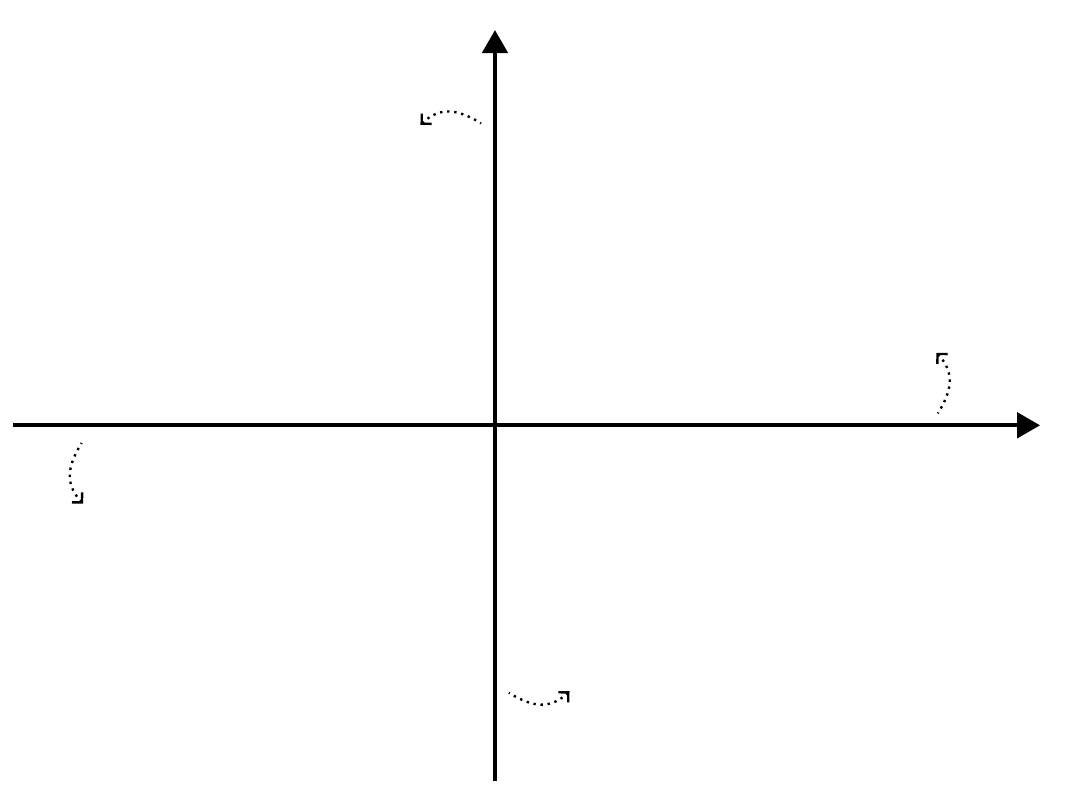}
			\put(97,30){$a_2$}
            \put(48,70){$b_2$}
			\put(75,45){Case $1.9$}
            \put(55,55){Case $1.10$}
            \put(25,55){Case $1.11$}
            \put(5,45){Case $1.12$}
            \put(5,20){Case $1.13$}
            \put(25,10){Case $1.14$}
            \put(55,10){Case $1.15$}
            \put(75,20){Case $1.16$}
		\end{overpic}
	\end{center}
	\caption{Bifurcation diagram in the variables $a_2$ and $b_2$, for $a_1 <0$ and $b_1=0$.}\label{figura_2.2}
\end{figure}

\noindent Here we will analyze the dynamics for all cases where $a_1 <0$ and $b_1=0$. See Figure \ref{figura_1.2} for the phase portraits corresponding to these cases and Figure \ref{figura_2.2} for the bifurcation diagram in the variables $a_2$ and $b_2$.

\vspace{0.5cm}

$\bullet$ \textbf{Case 1.9: $a_1<0, b_1 = 0, a_2 > 0, b_2 = 0$.} Given an arbitrary initial condition $(x_0,y_0)$, the trajectory has $\omega$-limit set $\{(-\infty,y_0)\}$ when $y_0 > \mu$ and has $\omega$-limit set $\{(+\infty,y_0)\}$ when $y_0 \leq \mu$.

$\bullet$ \textbf{Case 1.10: $a_1<0, b_1 = 0, a_2 > 0, b_2 > 0$.}  Given an arbitrary initial condition $(x_0,y_0)$,  the trajectory has $\omega$-limit set $\{(-\infty,y_0)\}$ when $y_0 > \mu$ and has $\omega$-limit set $\{(-\infty,\mu)\}$ when $y_0 \leq \mu$.

$\bullet$ \textbf{Case 1.11: $a_1<0, b_1 = 0, a_2=0, b_2 > 0$.} Given an arbitrary initial condition $(x_0,y_0)$, the trajectory has $\omega$-limit set $\{(-\infty,y_0)\}$ when $y_0 > \mu$ and has $\omega$-limit set $\{(-\infty,\mu)\}$ when $y_0 \leq \mu$.

$\bullet$ \textbf{Case 1.12: $a_1<0, b_1 = 0, a_2<0, b_2 > 0$.} Given an arbitrary initial condition $(x_0,y_0)$, the trajectory has $\omega$-limit set $\{(-\infty,y_0)\}$ when $y_0 > \mu$ and has $\omega$-limit set $\{(-\infty,\mu)\}$ when $y_0 \leq \mu$.

$\bullet$ \textbf{Case 1.13: $a_1<0, b_1 = 0, a_2 < 0, b_2 = 0$.} Given an arbitrary initial condition $(x_0,y_0)$, the trajectory has $\omega$-limit set $\{(-\infty,y_0)\}$.

$\bullet$ \textbf{Case 1.14: $a_1<0, b_1 = 0, a_2<0, b_2 < 0$.} Given an arbitrary initial condition $(x_0,y_0)$, the trajectory has $\omega$-limit set $\{(-\infty,y_0)\}$ when $y_0 > \mu$ and has $\omega$-limit set $\{(-\infty,-\infty)\}$ when $y_0 \leq \mu$.

$\bullet$ \textbf{Case 1.15: $a_1<0, b_1 = 0, a_2=0, b_2 < 0$.} Given an arbitrary initial condition $(x_0,y_0)$, the trajectory has $\omega$-limit set $\{(-\infty,y_0)\}$ when $y_0 > \mu$ and has $\omega$-limit set $\{(x_0,-\infty)\}$ when $y_0 \leq \mu$.

$\bullet$ \textbf{Case 1.16: $a_1<0, b_1 = 0, a_2>0, b_2 < 0$.} Given an arbitrary initial condition $(x_0,y_0)$, the trajectory has $\omega$-limit set $\{(-\infty,y_0)\}$ when $y_0 > \mu$ and has $\omega$-limit set $\{(+\infty,-\infty)\}$ when $y_0 \leq \mu$.

\vspace{0.5cm}

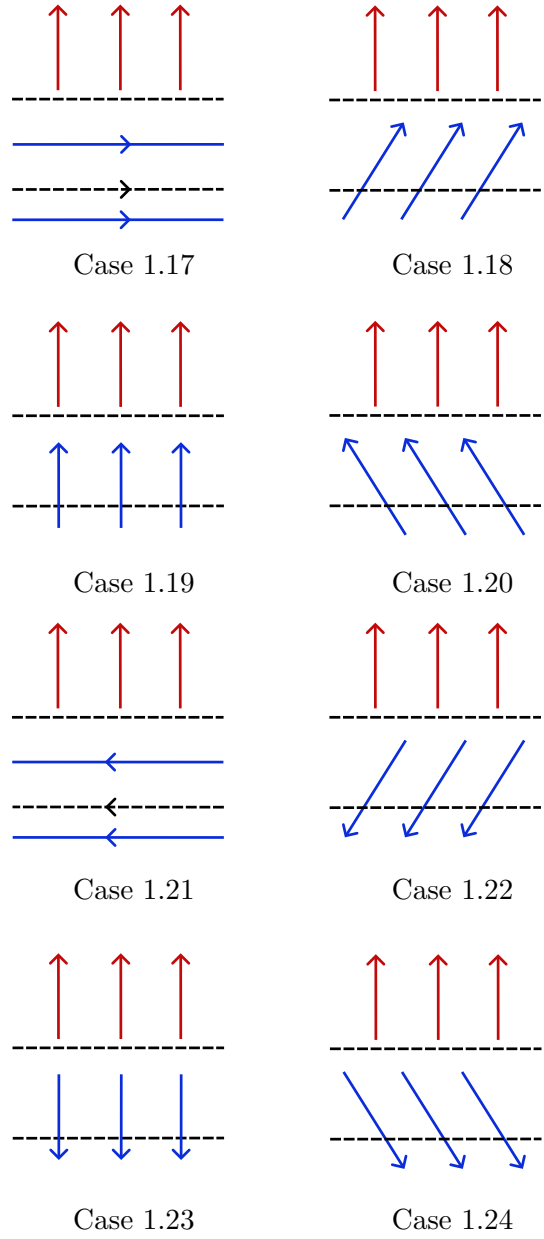
\begin{figure}[ht]
	\begin{center}
		\begin{overpic}[width=7cm]{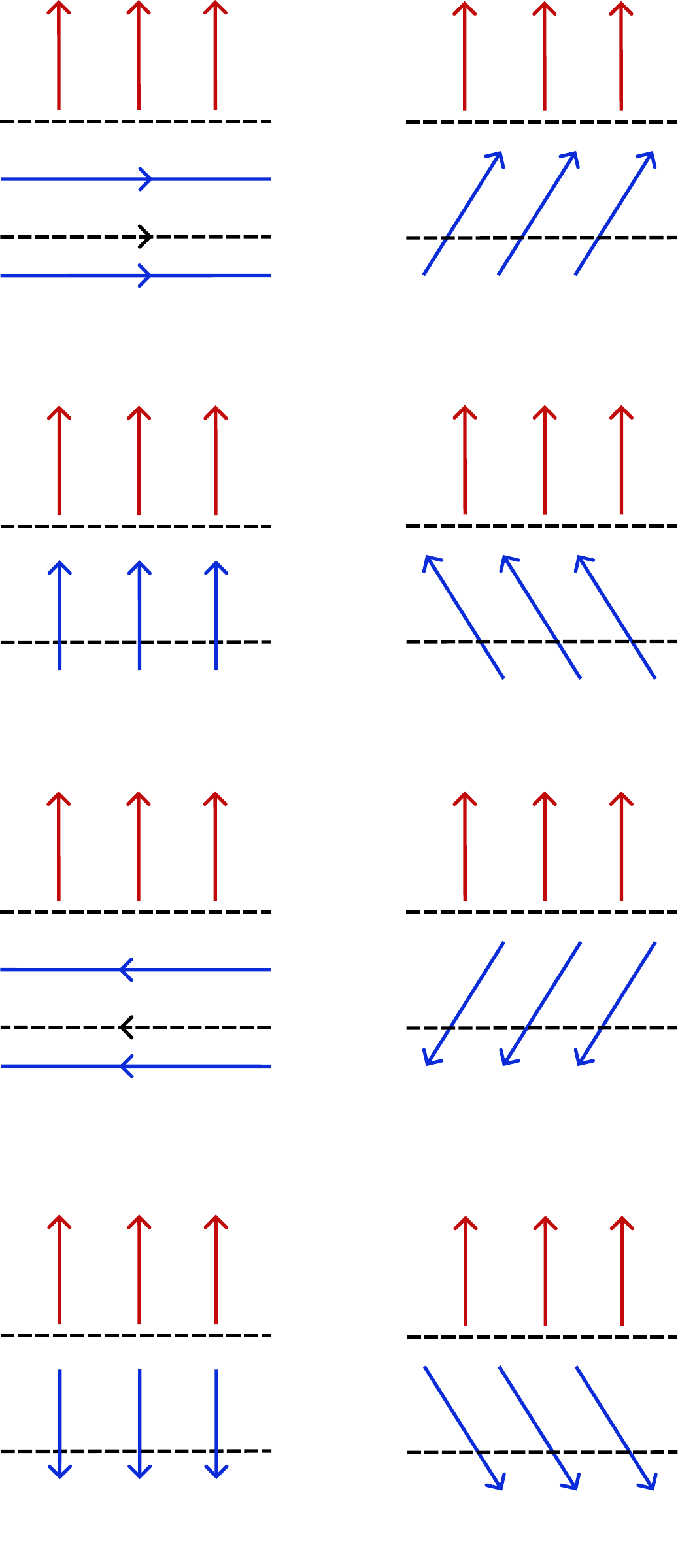}
			\put(5,78){Case $1.17$}
			\put(31,78){Case $1.18$}
			\put(5,52){Case $1.19$}
			\put(31,52){Case $1.20$}
			\put(5,27){Case $1.21$}
			\put(31,27){Case $1.22$}
            \put(5,0){Case $1.23$}
            \put(31,0){Case $1.24$}
		\end{overpic}
	\end{center}
	\caption{Dynamics for conditions $a_1=0$ and $b_1>0$.}\label{figura_1.3}
\end{figure}

\begin{figure}[ht]
	\begin{center}
		\begin{overpic}[width=9cm]{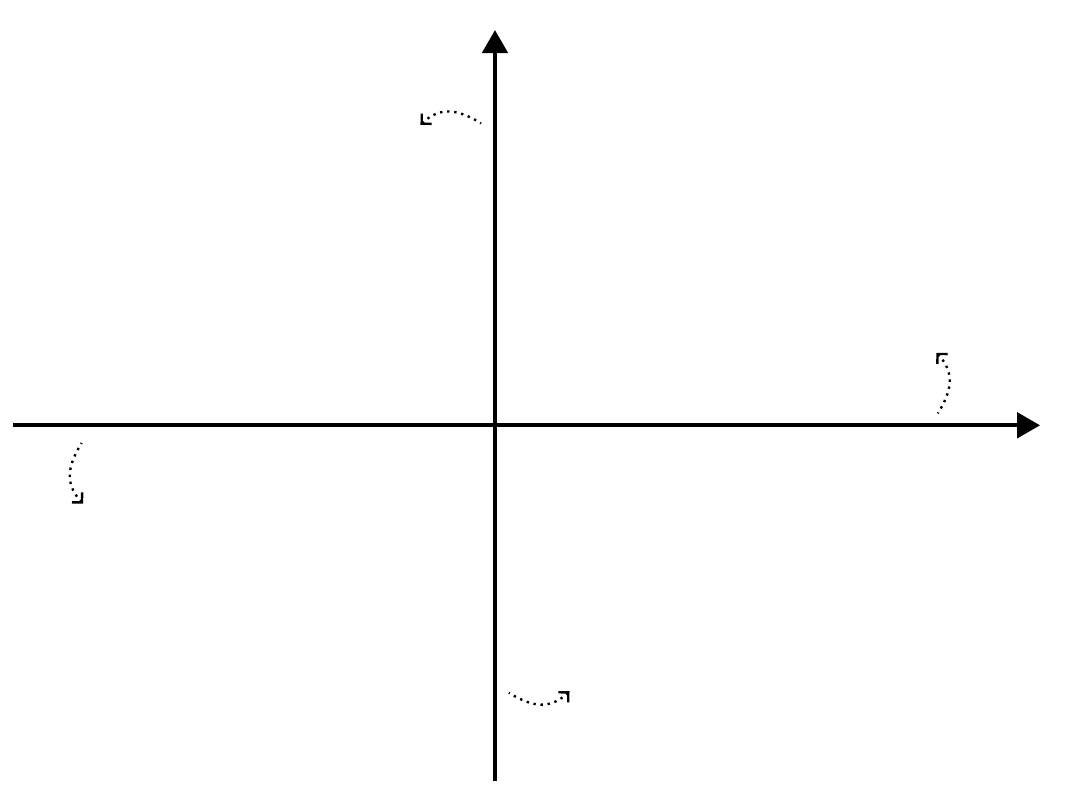}
			\put(97,30){$a_2$}
            \put(48,70){$b_2$}
			\put(75,45){Case $1.17$}
            \put(55,55){Case $1.18$}
            \put(25,55){Case $1.19$}
            \put(5,45){Case $1.20$}
            \put(5,20){Case $1.21$}
            \put(25,10){Case $1.22$}
            \put(55,10){Case $1.23$}
            \put(75,20){Case $1.24$}
		\end{overpic}
	\end{center}
	\caption{Bifurcation diagram in the variables $a_2$ and $b_2$, for $a_1 =0$ and $b_1>0$.}\label{figura_2.3}
\end{figure}

\noindent Here we will analyze the dynamics for all cases where $a_1 =0$ and $b_1>0$. See Figure \ref{figura_1.3} for the phase portraits corresponding to these cases and Figure \ref{figura_2.3} for the bifurcation diagram in the variables $a_2$ and $b_2$.

\vspace{0.5cm}

$\bullet$ \textbf{Case 1.17: $a_1=0, b_1 > 0, a_2 > 0, b_2 = 0$.} In this case, taking an arbitrary initial condition $(x_0,y_0)$, the trajectory has $\omega$-limit set $\{(x_0,+\infty)\}$ when $y_0 > \mu$ and has $\omega$-limit set $\{(+\infty,y_0)\}$ when $y_0 \leq \mu$.

$\bullet$ \textbf{Case 1.18: $a_1=0, b_1 > 0, a_2 > 0, b_2 > 0$.} In this case, taking an arbitrary initial condition $(x_0,y_0)$, the trajectory has $\omega$-limit set $\{(x_0,+\infty)\}$ when $y_0 > \mu$ and has $\omega$-limit set 
$$
\left\{\left(\dfrac{x_0 + a_2 (-y_0 + \mu)}{b_2},+\infty\right)\right\},
$$
when $y_0 \leq \mu$.

$\bullet$ \textbf{Case 1.19: $a_1=0, b_1 > 0, a_2 = 0, b_2 > 0$.} In this case, taking an arbitrary initial condition $(x_0,y_0)$, the trajectory has $\omega$-limit set $\{(x_0,+\infty)\}$.

$\bullet$ \textbf{Case 1.20: $a_1=0, b_1 > 0, a_2 < 0, b_2 > 0$.} In this case, taking an arbitrary initial condition $(x_0,y_0)$, the trajectory has $\omega$-limit set $\{(x_0,+\infty)\}$ when $y_0 > \mu$ and has $\omega$-limit set 
$$
\left\{\left(\dfrac{x_0 + a_2 (-y_0 + \mu)}{b_2},+\infty\right)\right\},
$$
when $y_0 \leq \mu$

$\bullet$ \textbf{Case 1.21: $a_1=0, b_1 > 0, a_2 < 0, b_2 = 0$.} In this case, taking an arbitrary initial condition $(x_0,y_0)$, the trajectory has $\omega$-limit set $\{(x_0,+\infty)\}$ when $y_0 > \mu$ and has $\omega$-limit set $\{(-\infty,y_0)\}$ when $y_0 \leq \mu$.

$\bullet$ \textbf{Case 1.22: $a_1=0, b_1 > 0, a_2 < 0, b_2 < 0$.} In this case, taking an arbitrary initial condition $(x_0,y_0)$, the trajectory has $\omega$-limit set $\{(x_0,+\infty)\}$ when $y_0 > \mu$ and has $\omega$-limit set $\{(-\infty,-\infty)\}$ when $y_0 \leq \mu$.

$\bullet$ \textbf{Case 1.23: $a_1=0, b_1 > 0, a_2 = 0, b_2 < 0$.} In this case, taking an arbitrary initial condition $(x_0,y_0)$, the trajectory has $\omega$-limit set $\{(x_0,+\infty)\}$ when $y_0 > \mu$ and has $\omega$-limit set $\{(x_0,-\infty)\}$ when $y_0 \leq \mu$.

$\bullet$ \textbf{Case 1.24: $a_1=0, b_1 > 0, a_2 > 0, b_2 < 0$.} In this case, taking an arbitrary initial condition $(x_0,y_0)$, the trajectory has $\omega$-limit set $\{(x_0,+\infty)\}$ when $y_0 > \mu$ and has $\omega$-limit set $\{(+\infty,-\infty)\}$ when $y_0 \leq \mu$.

\vspace{0.5cm}

\begin{figure}[ht]
	\begin{center}
		\begin{overpic}[width=7cm]{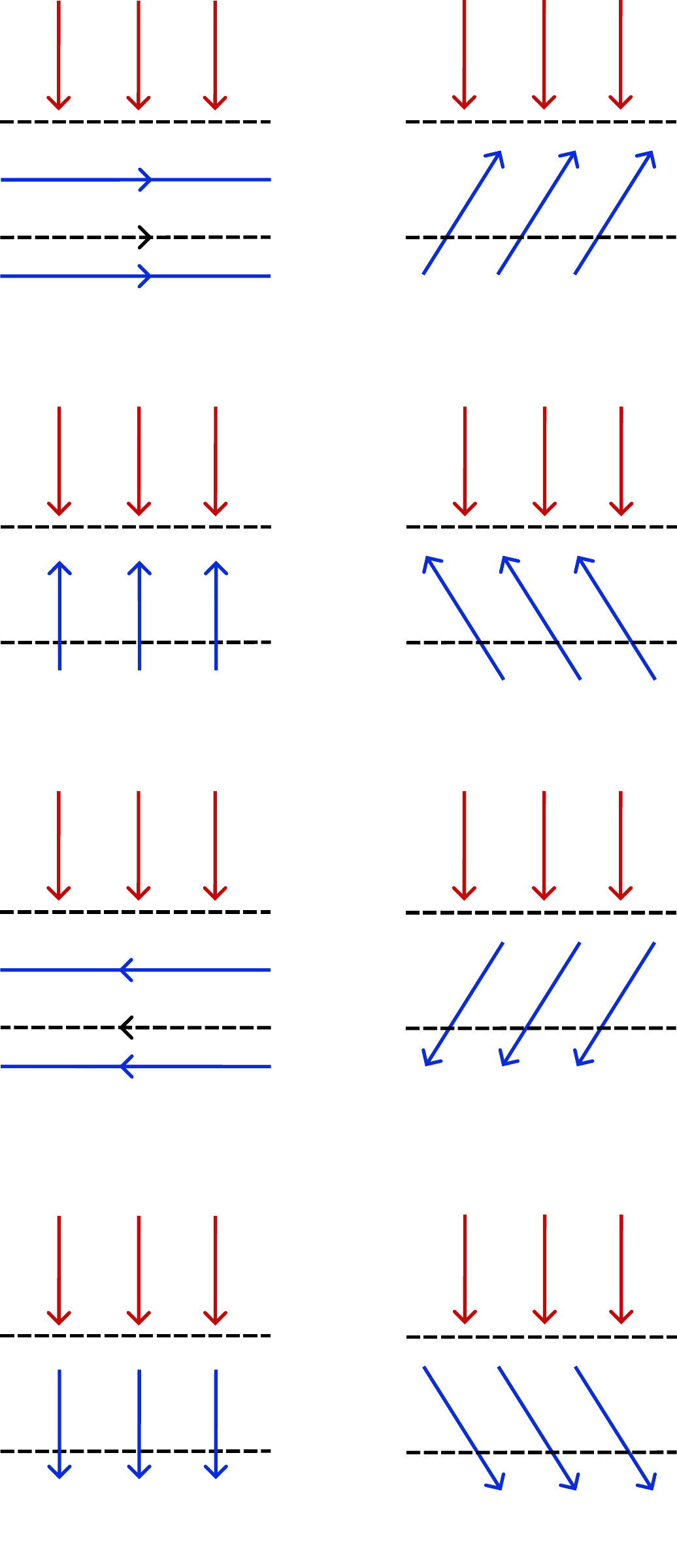}
			\put(5,78){Case $1.25$}
			\put(31,78){Case $1.26$}
			\put(5,52){Case $1.27$}
			\put(31,52){Case $1.28$}
			\put(5,27){Case $1.29$}
			\put(31,27){Case $1.30$}
            \put(5,0){Case $1.31$}
            \put(31,0){Case $1.32$}
		\end{overpic}
	\end{center}
	\caption{Dynamics for conditions $a_1=0$ and $b_1<0$.}\label{figura_1.4}
\end{figure}

\begin{figure}[ht]
	\begin{center}
		\begin{overpic}[width=9cm]{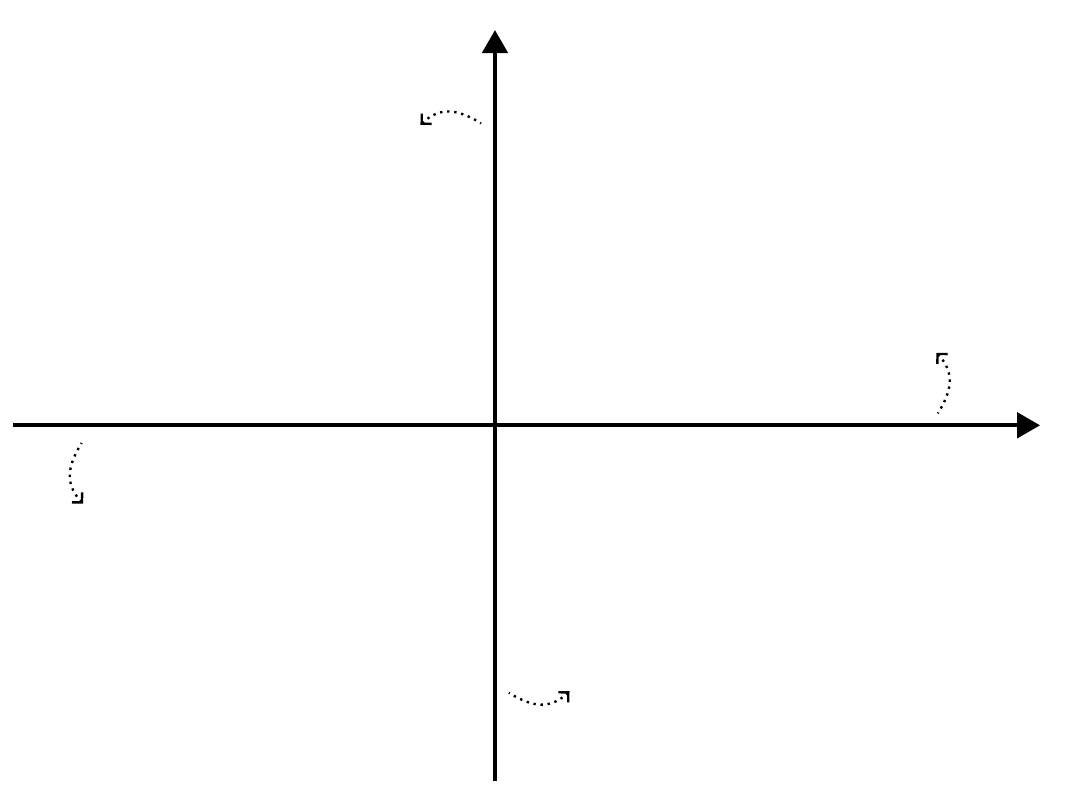}
			\put(97,30){$a_2$}
            \put(48,70){$b_2$}
			\put(75,45){Case $1.25$}
            \put(55,55){Case $1.26$}
            \put(25,55){Case $1.27$}
            \put(5,45){Case $1.28$}
            \put(5,20){Case $1.29$}
            \put(25,10){Case $1.30$}
            \put(55,10){Case $1.31$}
            \put(75,20){Case $1.32$}
		\end{overpic}
	\end{center}
	\caption{Bifurcation diagram in the variables $a_2$ and $b_2$, for $a_1 =0$ and $b_1<0$.}\label{figura_2.4}
\end{figure}

\noindent Here we will analyze the dynamics for all cases where $a_1 =0$ and $b_1<0$. See Figure \ref{figura_1.4} for the phase portraits corresponding to these cases and Figure \ref{figura_2.4} for the bifurcation diagram in the variables $a_2$ and $b_2$.

\vspace{0.5cm}

$\bullet$ \textbf{Case 1.25: $a_1=0, b_1 < 0, a_2 > 0, b_2 = 0$.} In this case, taking an arbitrary initial condition $(x_0,y_0)$, the trajectory has $\omega$-limit set $\{(+\infty,\mu)\}$ when $y_0 > \mu$ and has $\omega$-limit set $\{(+\infty,y_0)\}$ when $y_0 \leq \mu$.

$\bullet$ \textbf{Case 1.26: $a_1=0, b_1 < 0, a_2 > 0, b_2 > 0$.} To analyze the dynamics of the trajectories, we will use Poincaré first return map. See section \ref{subsec1}.

$\bullet$ \textbf{Case 1.27: $a_1=0, b_1 < 0, a_2 = 0, b_2 > 0$.} To analyze the dynamics of the trajectories, we will use Poincaré first return map. See section \ref{subsec1}.

$\bullet$ \textbf{Case 1.28: $a_1=0, b_1 < 0, a_2 < 0, b_2 > 0$.} To analyze the dynamics of the trajectories, we will use Poincaré first return map. See section \ref{subsec1}.

$\bullet$ \textbf{Case 1.29: $a_1=0, b_1 < 0, a_2 < 0, b_2 = 0$.} In this case, taking an arbitrary initial condition $(x_0,y_0)$, the trajectory has $\omega$-limit set $\{(-\infty,-\mu)\}$ when $y_0 > \mu$ and has $\omega$-limit set $\{(-\infty,y_0)\}$ when $y_0 \leq \mu$.

$\bullet$ \textbf{Case 1.30: $a_1=0, b_1 < 0, a_2 < 0, b_2 < 0$.} In this case, taking an arbitrary initial condition $(x_0,y_0)$, the trajectory has $\omega$-limit set $\{(-\infty,-\infty)\}$.

$\bullet$ \textbf{Case 1.31: $a_1=0, b_1 < 0, a_2 = 0, b_2 < 0$.} In this case, taking an arbitrary initial condition $(x_0,y_0)$, the trajectory has $\omega$-limit set $\{(x_0,-\infty)\}$.

$\bullet$ \textbf{Case 1.32: $a_1=0, b_1 < 0, a_2 > 0, b_2 < 0$.} In this case, taking an arbitrary initial condition $(x_0,y_0)$, the trajectory has $\omega$-limit set $\{(+\infty,-\infty)\}$.

\vspace{0.5cm}

\begin{figure}[ht]
	\begin{center}
		\begin{overpic}[width=7cm]{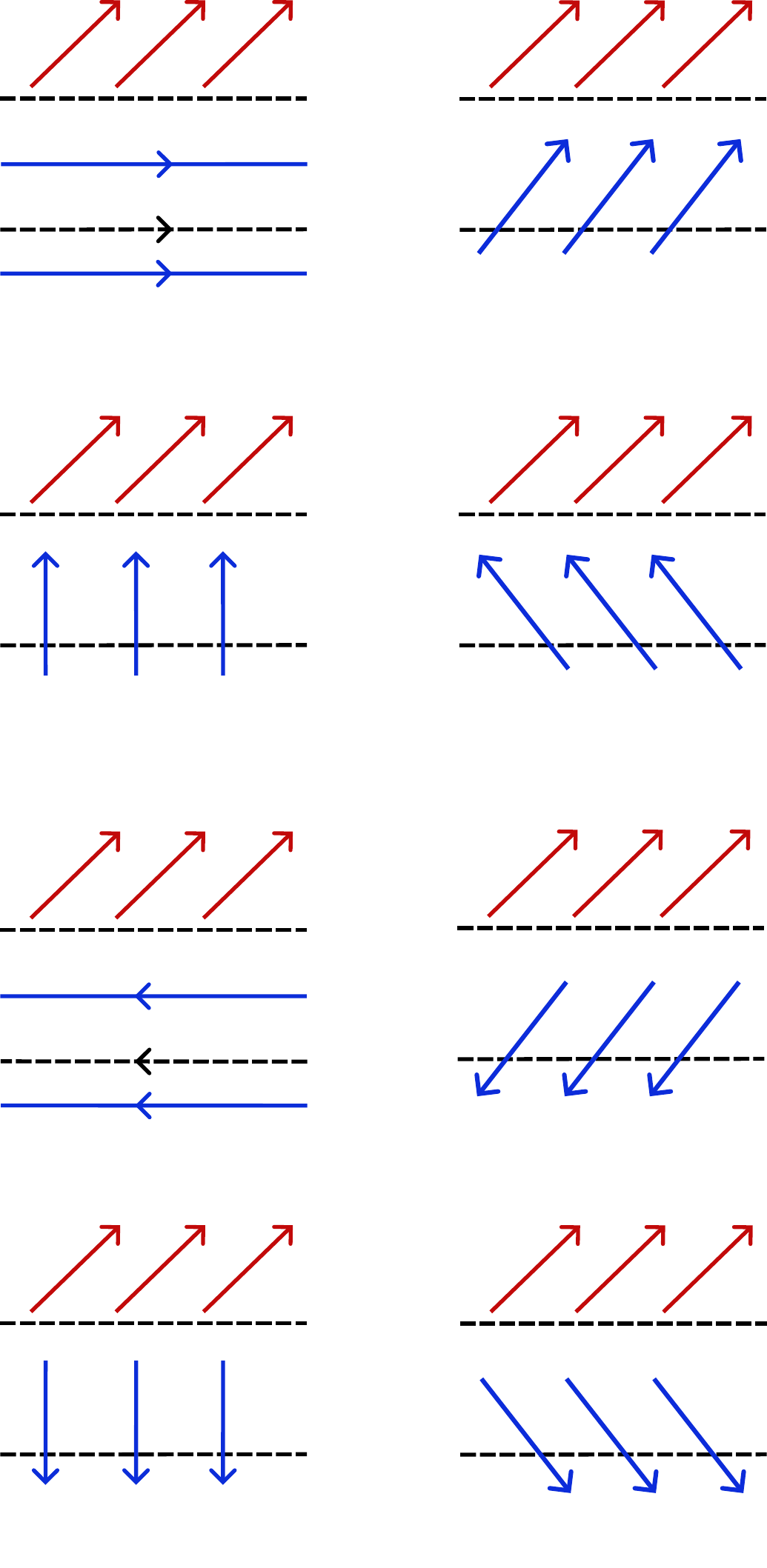}
			\put(5,78){Case $1.33$}
			\put(35,78){Case $1.34$}
			\put(5,52){Case $1.35$}
			\put(35,52){Case $1.36$}
			\put(5,25){Case $1.37$}
			\put(35,25){Case $1.38$}
            \put(5,0){Case $1.39$}
            \put(35,0){Case $1.40$}
		\end{overpic}
	\end{center}
	\caption{Dynamics for conditions $a_1>0$ and $b_1>0$.}\label{figura_1.5}
\end{figure}

\begin{figure}[ht]
	\begin{center}
		\begin{overpic}[width=9cm]{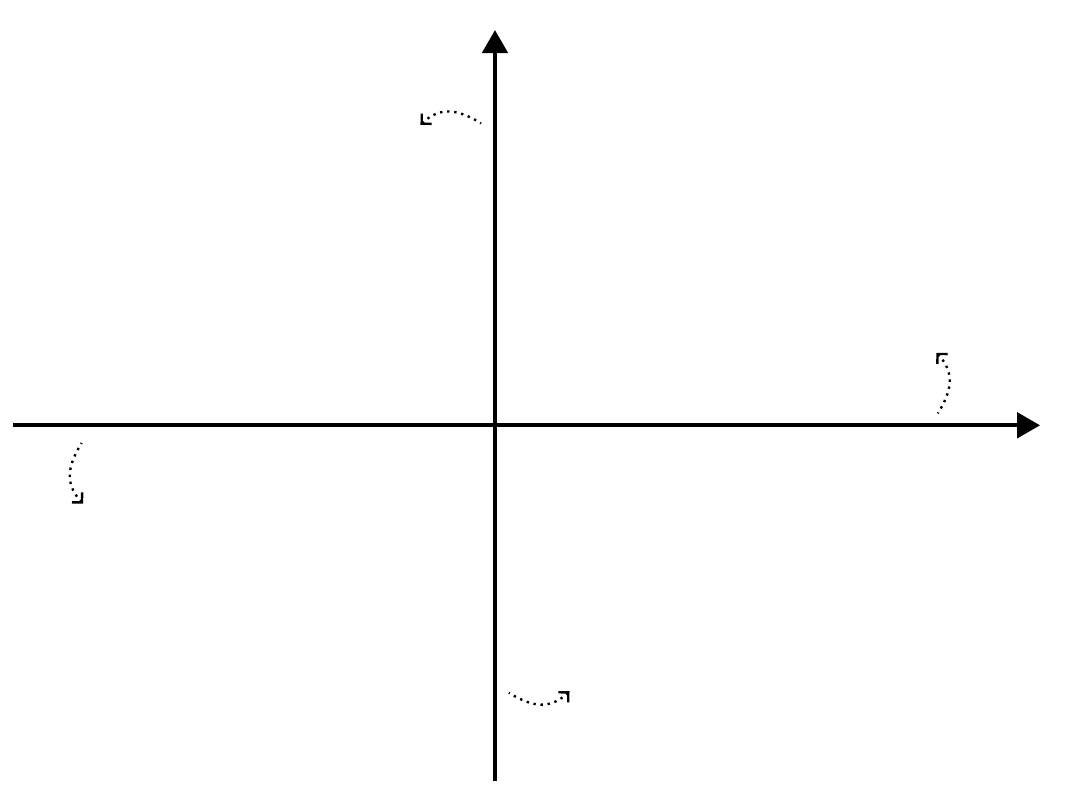}
			\put(97,30){$a_2$}
            \put(48,70){$b_2$}
			\put(75,45){Case $1.33$}
            \put(55,55){Case $1.34$}
            \put(25,55){Case $1.35$}
            \put(5,45){Case $1.36$}
            \put(5,20){Case $1.37$}
            \put(25,10){Case $1.38$}
            \put(55,10){Case $1.39$}
            \put(75,20){Case $1.40$}
		\end{overpic}
	\end{center}
	\caption{Bifurcation diagram in the variables $a_2$ and $b_2$, for $a_1 >0$ and $b_1>0$.}\label{figura_2.5}
\end{figure}

\noindent Here we will analyze the dynamics for all cases where $a_1 >0$ and $b_1>0$. See Figure \ref{figura_1.5} for the phase portraits corresponding to these cases and Figure \ref{figura_2.5} for the bifurcation diagram in the variables $a_2$ and $b_2$.

\vspace{0.5cm}

$\bullet$ \textbf{Case 1.33: $a_1>0, b_1 > 0, a_2 > 0, b_2 = 0$.} In this case, taking an arbitrary initial condition $(x_0,y_0)$, the trajectory has $\omega$-limit set $\{(+\infty,+\infty)\}$ when $y_0 > \mu$ and has $\omega$-limit set $\{(+\infty,y_0)\}$ when $y_0 \leq \mu$.

$\bullet$ \textbf{Case 1.34: $a_1>0, b_1 > 0, a_2>0, b_2 > 0$.} In this case, taking an arbitrary initial condition $(x_0,y_0)$, the trajectory has $\omega$-limit set $\{(+\infty,+\infty)\}$.

$\bullet$ \textbf{Case 1.35: $a_1>0, b_1 > 0, a_2=0, b_2 > 0$.} In this case, taking an arbitrary initial condition $(x_0,y_0)$, the trajectory has $\omega$-limit set $\{(+\infty,+\infty)\}$.

$\bullet$ \textbf{Case 1.36: $a_1>0, b_1 > 0, a_2<0, b_2 > 0$.} In this case, taking an arbitrary initial condition $(x_0,y_0)$, the trajectory has $\omega$-limit set $\{(+\infty,+\infty)\}$.

$\bullet$ \textbf{Case 1.37: $a_1>0, b_1 > 0, a_2 < 0, b_2 = 0$.} In this case, taking an arbitrary initial condition $(x_0,y_0)$, the trajectory has $\omega$-limit set $\{(+\infty,+\infty)\}$ when $y_0 > \mu$ and has $\omega$-limit set $\{(-\infty,y_0)\}$ when $y_0 \leq \mu$.

$\bullet$ \textbf{Case 1.38: $a_1>0, b_1 > 0, a_2<0, b_2 < 0$.} In this case, taking an arbitrary initial condition $(x_0,y_0)$, the trajectory has $\omega$-limit set $\{(+\infty,+\infty)\}$ when $y_0 > \mu$ and has $\omega$-limit set $\{(-\infty,-\infty)\}$ when $y_0 \leq \mu$.

$\bullet$ \textbf{Case 1.39: $a_1>0, b_1 > 0, a_2 = 0, b_2 < 0$.} In this case, taking an arbitrary initial condition $(x_0,y_0)$, the trajectory has $\omega$-limit set $\{(+\infty,+\infty)\}$ when $y_0 > \mu$ and has $\omega$-limit set $\{(x_0,-\infty)\}$ when $y_0 \leq \mu$.

$\bullet$ \textbf{Case 1.40: $a_1>0, b_1 > 0, a_2>0, b_2 < 0$.} In this case, taking an arbitrary initial condition $(x_0,y_0)$, the trajectory has $\omega$-limit set $\{(+\infty,+\infty)\}$ when $y_0 > \mu$ and has $\omega$-limit set $\{(+\infty,-\infty)\}$ when $y_0 \leq \mu$.

\vspace{0.5cm}

\begin{figure}[ht]
	\begin{center}
		\begin{overpic}[width=7cm]{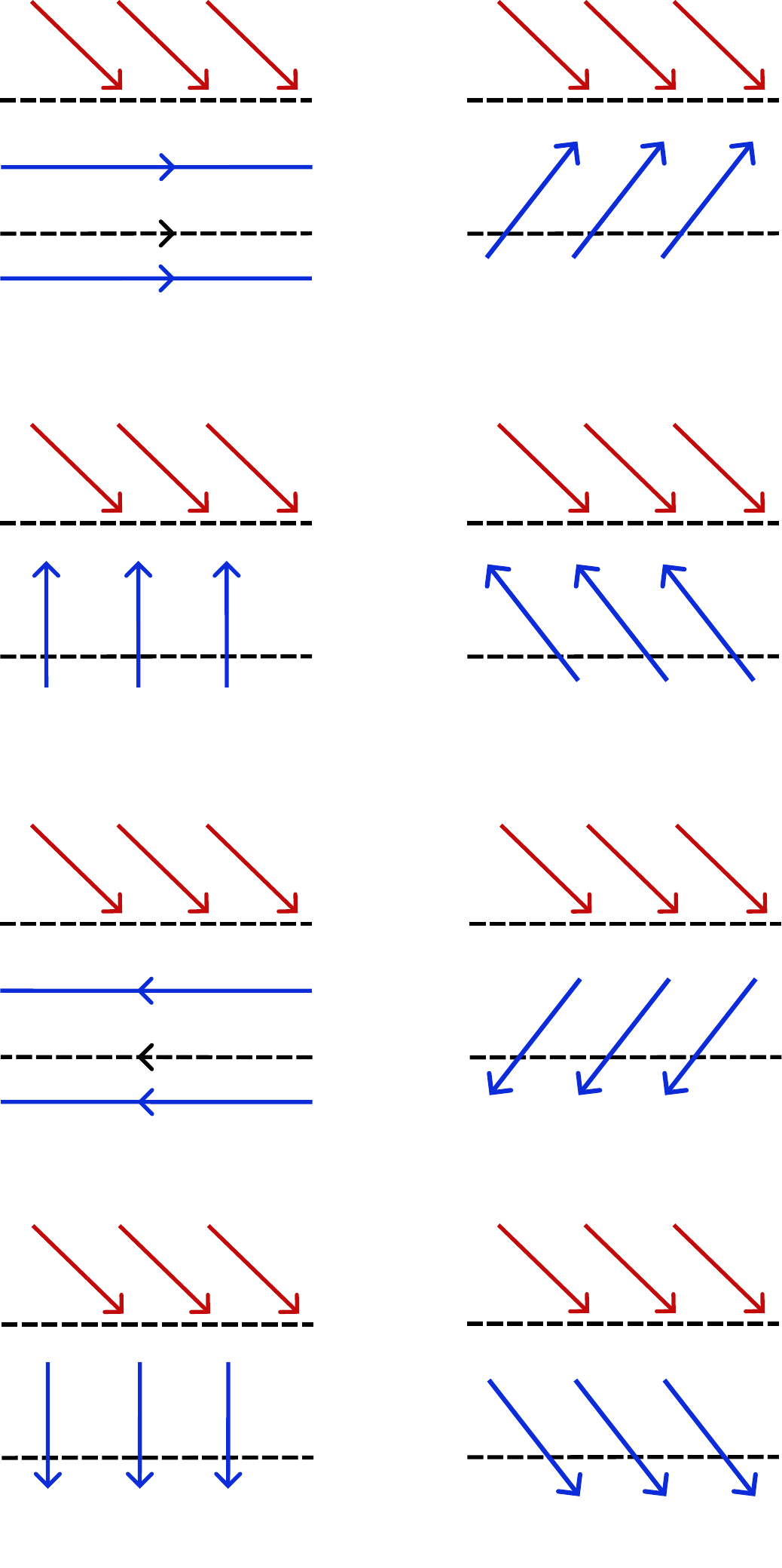}
			\put(5,78){Case $1.41$}
			\put(35,78){Case $1.42$}
			\put(5,52){Case $1.43$}
			\put(35,52){Case $1.44$}
			\put(5,25){Case $1.45$}
			\put(35,25){Case $1.46$}
            \put(5,0){Case $1.47$}
            \put(35,0){Case $1.48$}
		\end{overpic}
	\end{center}
	\caption{Dynamics for conditions $a_1>0$ and $b_1<0$.}\label{figura_1.6}
\end{figure}

\begin{figure}[ht]
	\begin{center}
		\begin{overpic}[width=9cm]{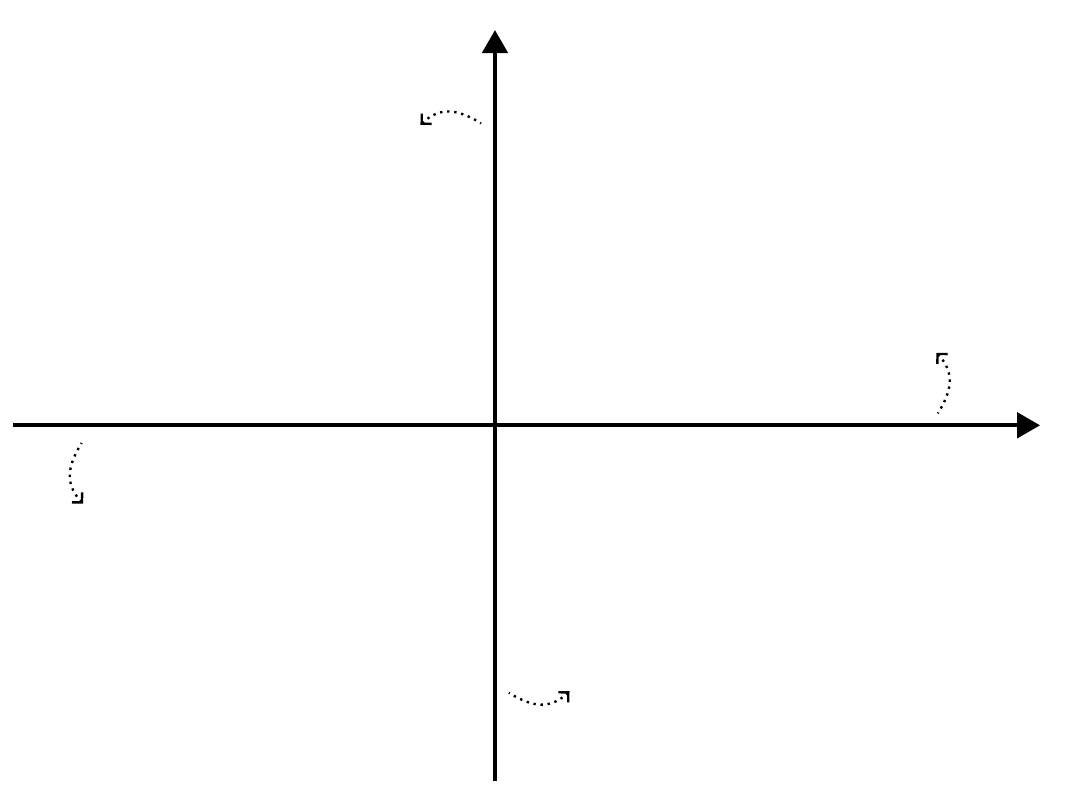}
			\put(97,30){$a_2$}
            \put(48,70){$b_2$}
			\put(75,45){Case $1.41$}
            \put(55,55){Case $1.42$}
            \put(25,55){Case $1.43$}
            \put(5,45){Case $1.44$}
            \put(5,20){Case $1.45$}
            \put(25,10){Case $1.46$}
            \put(55,10){Case $1.47$}
            \put(75,20){Case $1.48$}
		\end{overpic}
	\end{center}
	\caption{Bifurcation diagram in the variables $a_2$ and $b_2$, for $a_1 >0$ and $b_1<0$.}\label{figura_2.6}
\end{figure}

\noindent Here we will analyze the dynamics for all cases where $a_1 >0$ and $b_1<0$. See Figure \ref{figura_1.6} for the phase portraits corresponding to these cases and Figure \ref{figura_2.6} for the bifurcation diagram in the variables $a_2$ and $b_2$.

\vspace{0.5cm}

$\bullet$ \textbf{Case 1.41: $a_1>0, b_1 < 0, a_2 > 0, b_2 = 0$.} In this case, taking an arbitrary initial condition $(x_0,y_0)$, the trajectory has $\omega$-limit set $\{(+\infty,-\mu)\}$ when $y_0 > \mu$ and has $\omega$-limit set $\{(+\infty,y_0)\}$ when $y_0 \leq \mu$.

$\bullet$ \textbf{Case 1.42: $a_1>0, b_1 < 0, a_2 > 0, b_2 > 0$.} To analyze the dynamics of the trajectories, we will use Poincaré first return map. See section \ref{subsec1}.

$\bullet$ \textbf{Case 1.43: $a_1>0, b_1 < 0, a_2=0, b_2 > 0$.} To analyze the dynamics of the trajectories, we will use Poincaré first return map. See section \ref{subsec1}.

$\bullet$ \textbf{Case 1.44: $a_1>0, b_1 < 0, a_2 < 0, b_2 > 0$.} To analyze the dynamics of the trajectories, we will use Poincaré first return map. See section \ref{subsec1}.

$\bullet$ \textbf{Case 1.45: $a_1>0, b_1 < 0, a_2 < 0, b_2 = 0$.} In this case, taking an arbitrary initial condition $(x_0,y_0)$, the trajectory has $\omega$-limit set $\{(-\infty,-\mu)\}$ when $y_0 > \mu$ and has $\omega$-limit set $\{(-\infty,y_0)\}$ when $y_0 \leq \mu$.

$\bullet$ \textbf{Case 1.46: $a_1>0, b_1 < 0, a_2 < 0, b_2 < 0$.} In this case, taking an arbitrary initial condition $(x_0,y_0)$, the trajectory has $\omega$-limit set $\{(-\infty,-\infty)\}$.

$\bullet$ \textbf{Case 1.47: $a_1>0, b_1 < 0, a_2 = 0, b_2 < 0$.} In this case, taking an arbitrary initial condition $(x_0,y_0)$, the trajectory has $\omega$-limit set  
$$
\left\{\left(\dfrac{x_0 + a_1 (-y_0 + \mu)}{b_1},-\infty\right)\right\},
$$
when $y_0 > \mu$ and has $\omega$-limit set $\{(x_0,-\infty)\}$ when $y_0 \leq \mu$.

$\bullet$ \textbf{Case 1.48: $a_1>0, b_1 < 0, a_2 > 0, b_2 < 0$.} In this case, taking an arbitrary initial condition $(x_0,y_0)$, the trajectory has $\omega$-limit set $\{(+\infty,-\infty)\}$.

\vspace{0.5cm}

\begin{figure}[ht]
	\begin{center}
		\begin{overpic}[width=7cm]{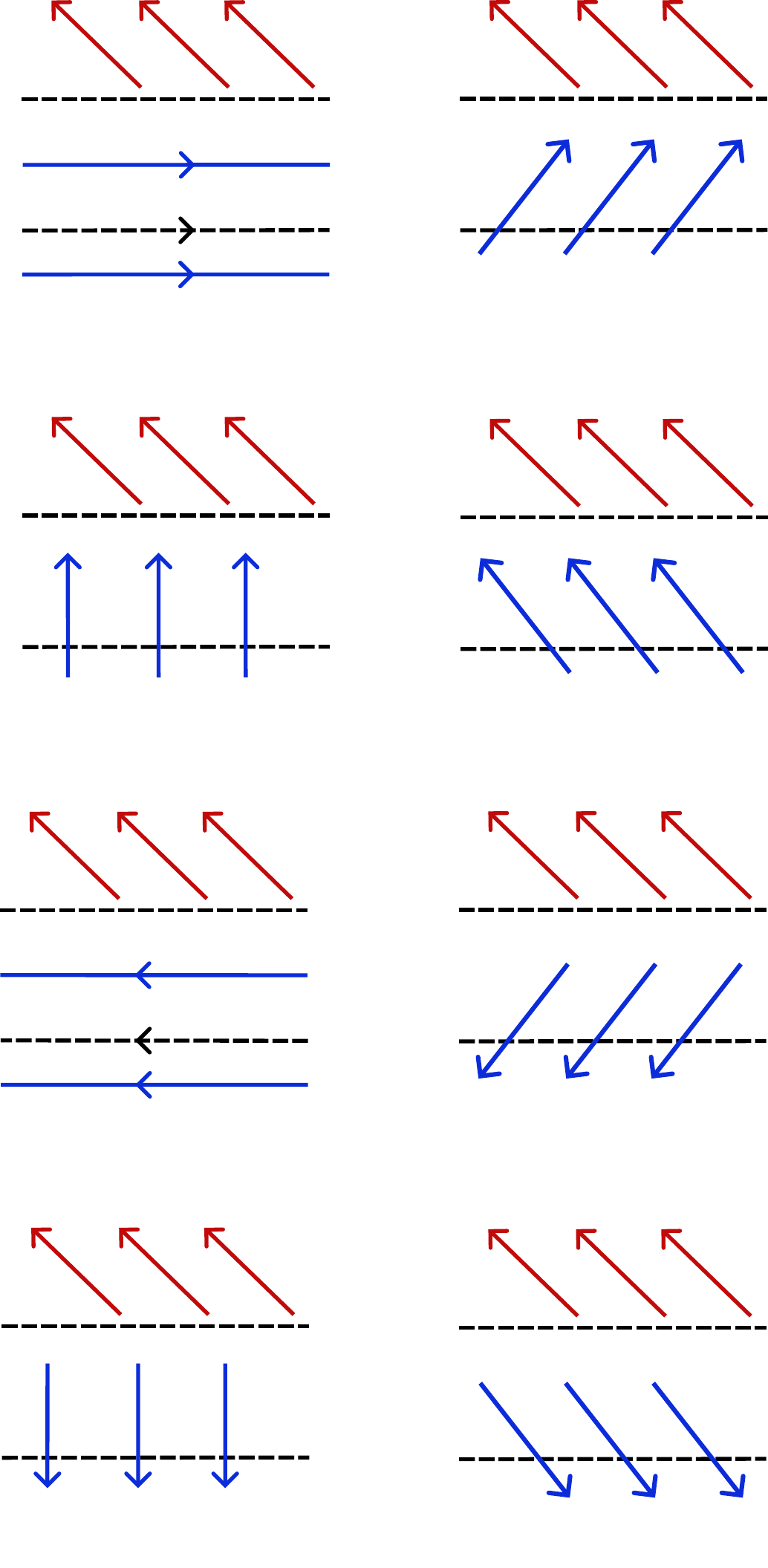}
			\put(5,78){Case $1.49$}
			\put(35,78){Case $1.50$}
			\put(5,52){Case $1.51$}
			\put(35,52){Case $1.52$}
			\put(5,25){Case $1.53$}
			\put(35,25){Case $1.54$}
            \put(5,0){Case $1.55$}
            \put(35,0){Case $1.56$}
		\end{overpic}
	\end{center}
	\caption{Dynamics for conditions $a_1<0$ and $b_1>0$.}\label{figura_1.7}
\end{figure}

\begin{figure}[ht]
	\begin{center}
		\begin{overpic}[width=9cm]{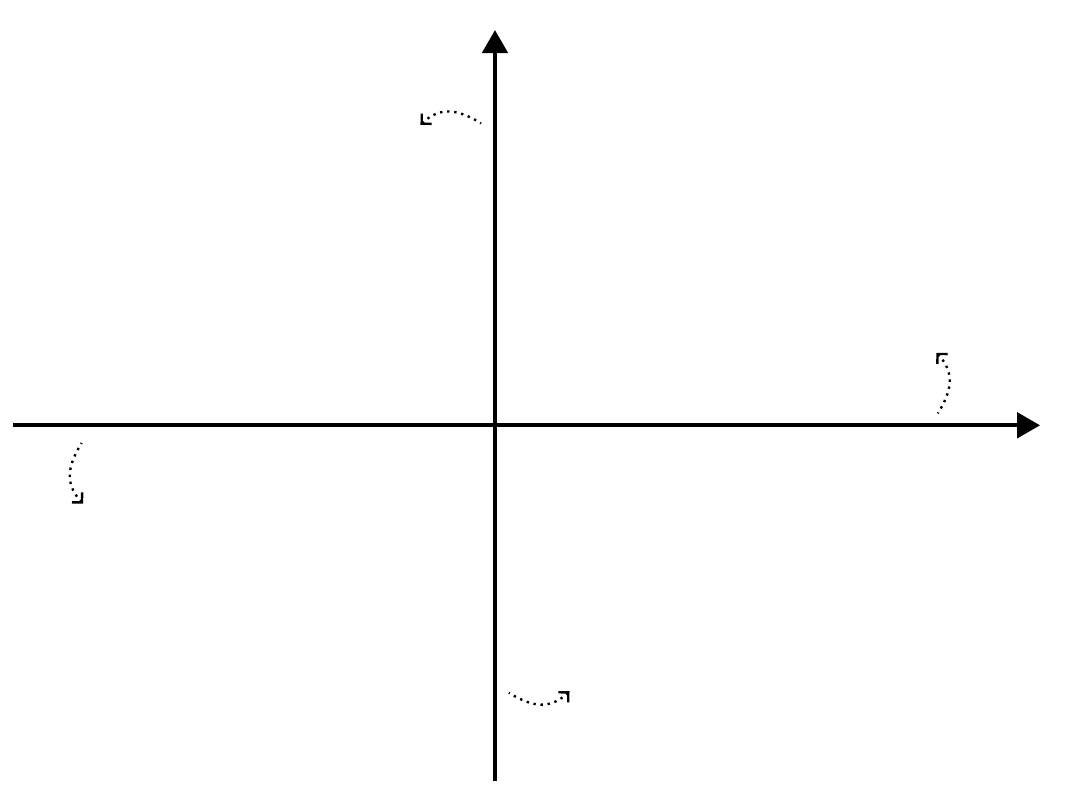}
			\put(97,30){$a_2$}
            \put(48,70){$b_2$}
			\put(75,45){Case $1.49$}
            \put(55,55){Case $1.50$}
            \put(25,55){Case $1.51$}
            \put(5,45){Case $1.52$}
            \put(5,20){Case $1.53$}
            \put(25,10){Case $1.54$}
            \put(55,10){Case $1.55$}
            \put(75,20){Case $1.56$}
		\end{overpic}
	\end{center}
	\caption{Bifurcation diagram in the variables $a_2$ and $b_2$, for $a_1 <0$ and $b_1> 0$.}\label{figura_2.7}
\end{figure}

\noindent Here we will analyze the dynamics for all cases where $a_1 <0$ and $b_1>0$. See Figure \ref{figura_1.7} for the phase portraits corresponding to these cases and Figure \ref{figura_2.7} for the bifurcation diagram in the variables $a_2$ and $b_2$.

\vspace{0.5cm}

$\bullet$ \textbf{Case 1.49: $a_1<0, b_1 > 0, a_2 > 0, b_2 = 0$.} In this case, taking an arbitrary initial condition $(x_0,y_0)$, the trajectory has $\omega$-limit set $\{(-\infty,+\infty)\}$ when $y_0 > \mu$ and has $\omega$-limit set $\{(+\infty,y_0)\}$ when $y_0 \leq \mu$.

$\bullet$ \textbf{Case 1.50: $a_1<0, b_1 > 0, a_2 > 0, b_2 > 0$.} In this case, taking an arbitrary initial condition $(x_0,y_0)$, the trajectory has $\omega$-limit set $\{(-\infty,+\infty)\}$.

$\bullet$ \textbf{Case 1.51: $a_1<0, b_1 > 0, a_2=0, b_2 > 0$.} In this case, taking an arbitrary initial condition $(x_0,y_0)$, the trajectory has $\omega$-limit set $\{(-\infty,+\infty)\}$.

$\bullet$ \textbf{Case 1.52: $a_1<0, b_1 > 0, a_2 < 0, b_2 > 0$.} In this case, taking an arbitrary initial condition $(x_0,y_0)$, the trajectory has $\omega$-limit set $\{(-\infty,+\infty)\}$.

$\bullet$ \textbf{Case 1.53: $a_1<0, b_1 > 0, a_2 < 0, b_2 = 0$.} In this case, taking an arbitrary initial condition $(x_0,y_0)$, the trajectory has $\omega$-limit set $\{(-\infty,+\infty)\}$ when $y_0 > \mu$ and has $\omega$-limit set $\{(-\infty,y_0)\}$ when $y_0 \leq \mu$.

$\bullet$ \textbf{Case 1.54: $a_1<0, b_1 > 0, a_2 < 0, b_2 < 0$.} In this case, taking an arbitrary initial condition $(x_0,y_0)$, the trajectory has $\omega$-limit set $\{(-\infty,+\infty)\}$ when $y_0 > \mu$ and has $\omega$-limit set $\{(-\infty,-\infty)\}$ when $y_0 \leq \mu$.

$\bullet$ \textbf{Case 1.55: $a_1<0, b_1 > 0, a_2 = 0, b_2 < 0$.} In this case, taking an arbitrary initial condition $(x_0,y_0)$, the trajectory has $\omega$-limit set $\{(-\infty,+\infty)\}$ when $y_0 > \mu$ and has $\omega$-limit set $\{(x_0,-\infty)\}$ when $y_0 \leq \mu$.

$\bullet$ \textbf{Case 1.56: $a_1<0, b_1 > 0, a_2 > 0, b_2 < 0$.} In this case, taking an arbitrary initial condition $(x_0,y_0)$, the trajectory has $\omega$-limit set $\{(-\infty,+\infty)\}$ when $y_0 > \mu$ and has $\omega$-limit set $\{(+\infty,-\infty)\}$ when $y_0 \leq \mu$.

\vspace{0.5cm}

\begin{figure}[ht]
	\begin{center}
		\begin{overpic}[width=7cm]{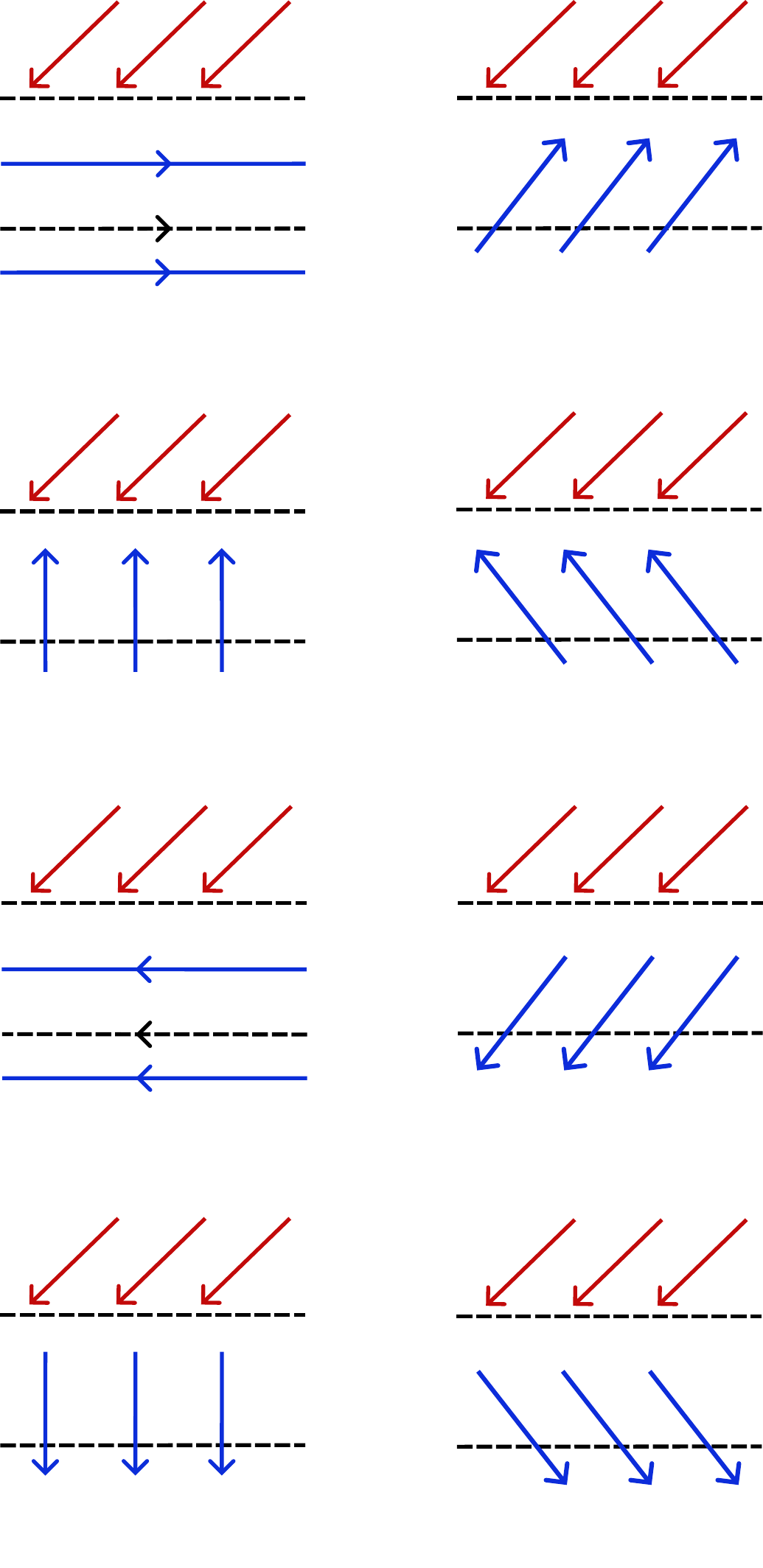}
			\put(5,78){Case $1.57$}
			\put(35,78){Case $1.58$}
			\put(5,52){Case $1.59$}
			\put(35,52){Case $1.60$}
			\put(5,25){Case $1.61$}
			\put(35,25){Case $1.62$}
            \put(5,0){Case $1.63$}
            \put(35,0){Case $1.64$}
		\end{overpic}
	\end{center}
	\caption{Dynamics for conditions $a_1<0$ and $b_1<0$.}\label{figura_1.8}
\end{figure}

\begin{figure}[ht]
	\begin{center}
		\begin{overpic}[width=9cm]{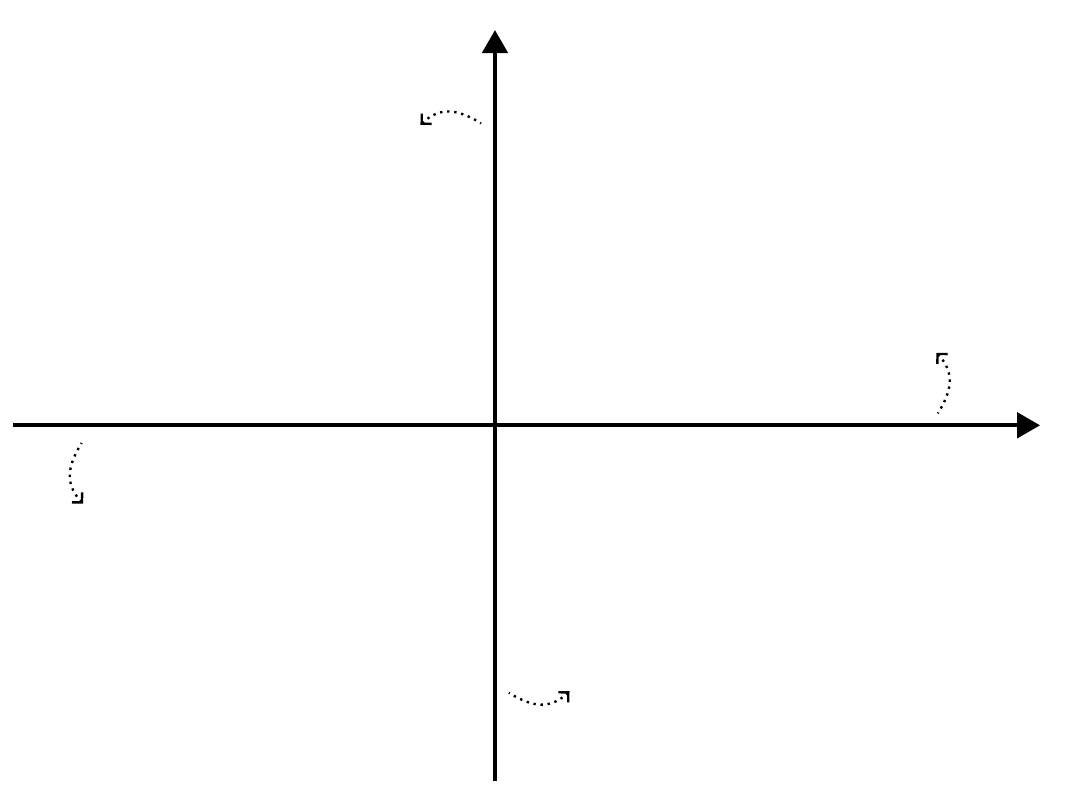}
			\put(97,30){$a_2$}
            \put(48,70){$b_2$}
			\put(75,45){Case $1.57$}
            \put(55,55){Case $1.58$}
            \put(25,55){Case $1.59$}
            \put(5,45){Case $1.60$}
            \put(5,20){Case $1.61$}
            \put(25,10){Case $1.62$}
            \put(55,10){Case $1.63$}
            \put(75,20){Case $1.64$}
		\end{overpic}
	\end{center}
	\caption{Bifurcation diagram in the variables $a_2$ and $b_2$, for $a_1 < 0$ and $b_1<0$.}\label{figura_2.8}
\end{figure}

\noindent Here we will analyze the dynamics for all cases where $a_1 <0$ and $b_1<0$. See Figure \ref{figura_1.8} for the phase portraits corresponding to these cases and Figure \ref{figura_2.8} for the bifurcation diagram in the variables $a_2$ and $b_2$.

\vspace{0.5cm}

$\bullet$ \textbf{Case 1.57: $a_1<0, b_1 < 0, a_2 > 0, b_2 = 0$.} In this case, taking an arbitrary initial condition $(x_0,y_0)$, the trajectory has $\omega$-limit set $\{(+\infty,\mu)\}$ when $y_0 > \mu$ and has $\omega$-limit set $\{(+\infty,y_0)\}$ when $y_0 \leq \mu$.

$\bullet$ \textbf{Case 1.58: $a_1<0, b_1 < 0, a_2>0, b_2 > 0$.} To analyze the dynamics of the trajectories, we will use Poincaré first return map. See section \ref{subsec1}.

$\bullet$ \textbf{Case 1.59: $a_1<0, b_1 < 0, a_2=0, b_2 > 0$.} To analyze the dynamics of the trajectories, we will use Poincaré first return map. See section \ref{subsec1}.

$\bullet$ \textbf{Case 1.60: $a_1<0, b_1 < 0, a_2<0, b_2 > 0$.} To analyze the dynamics of the trajectories, we will use Poincaré first return map. See section \ref{subsec1}.

$\bullet$ \textbf{Case 1.61: $a_1<0, b_1 < 0, a_2 < 0, b_2 = 0$.} In this case, taking an arbitrary initial condition $(x_0,y_0)$, the trajectory has $\omega$-limit set $\{(-\infty,\mu)\}$ when $y_0 > \mu$ and has $\omega$-limit set $\{(-\infty,y_0)\}$ when $y_0 \leq \mu$.

$\bullet$ \textbf{Case 1.62: $a_1<0, b_1 < 0, a_2<0, b_2 < 0$.} In this case, taking an arbitrary initial condition $(x_0,y_0)$, the trajectory has $\omega$-limit set $\{(-\infty,-\infty)\}$.

$\bullet$ \textbf{Case 1.63: $a_1<0, b_1 < 0, a_2 = 0, b_2 < 0$.} In this case, taking an arbitrary initial condition $(x_0,y_0)$, the trajectory has $\omega$-limit set 
$$
\left\{\left(\dfrac{x_0 + a_1 (-y_0 + \mu)}{b_1},-\infty\right)\right\},
$$
when $y_0 > \mu$ and has $\omega$-limit set $\{(x_0,-\infty)\}$ when $y_0 \leq \mu$.

$\bullet$ \textbf{Case 1.64: $a_1<0, b_1 < 0, a_2>0, b_2 < 0$.} In this case, taking an arbitrary initial condition $(x_0,y_0)$, the trajectory has $\omega$-limit set $\{(+\infty,-\infty)\}$.

\subsubsection{\textbf{Poincaré first return map}}\label{subsec1}

In this section, via the Poincaré first return map, let us describe the dynamics of the trajectories of $X_1$ and $X_2$ in the cases 1.26, 1.27, 1.28, 1.42, 1.43, 1.44, 1.58, 1.59, 1.60. Note that $b_1b_2\neq0$.

Let us suppose an initial point  $p_1=\left(x_0^1, \mu\right)$ located at the superior boundary of the hysteresis band and under the action of the vector field $X_1$. The solutions for this vector field are given by
$$
x_1\left(t\right)=a_1t+x_0^1 \quad \textrm{and} \quad y_1\left(t\right)=b_1t+\mu. 
$$
The time which the flow takes to reach the another boundary of the hysteresis band $y=-\mu$ is given by 
\begin{equation}\label{eq1}
	t=-\dfrac{2\mu}{b_1}.
\end{equation}
Substituting \eqref{eq1} in the component solution $x_1$ gives us
\[
x_1\left(-\dfrac{2\mu}{b_1}\right):=x_0^2=x_0^1-\dfrac{2a_1\mu}{b_1}.
\]
From this moment, the flow is governed by the vector field $X_2$ located under the hysteresis region and the solutions for this vector field from the initial position $p_2=(x_0^2,-\mu)$ are given by

\[
	x_2(t)=a_2t+x_0^1-\dfrac{2a_1\mu}{b_1} \quad \textrm{and} \quad y_2\left(t\right)=b_2t-\mu.
\]
The time that this flow takes to achieve the superior boundary $y=\mu$ is then given by $t=2\mu/b_2$ and substituting this time $t$ in $x_2$ gives
\[
x_2\left(\dfrac{2\mu}{b_2}\right)=x_0^1-\dfrac{2a_1\mu}{b_1}+\dfrac{2a_2\mu}{b_2}.
\]
Hence, the Poincaré first return map is given by 
\[
	P\left(x\right)=x+2\mu\left(\frac{a_2}{b_2}-\frac{a_1}{b_1}\right).
\]

\begin{figure}[ht]
	\begin{center}
		\begin{overpic}[width=4cm]{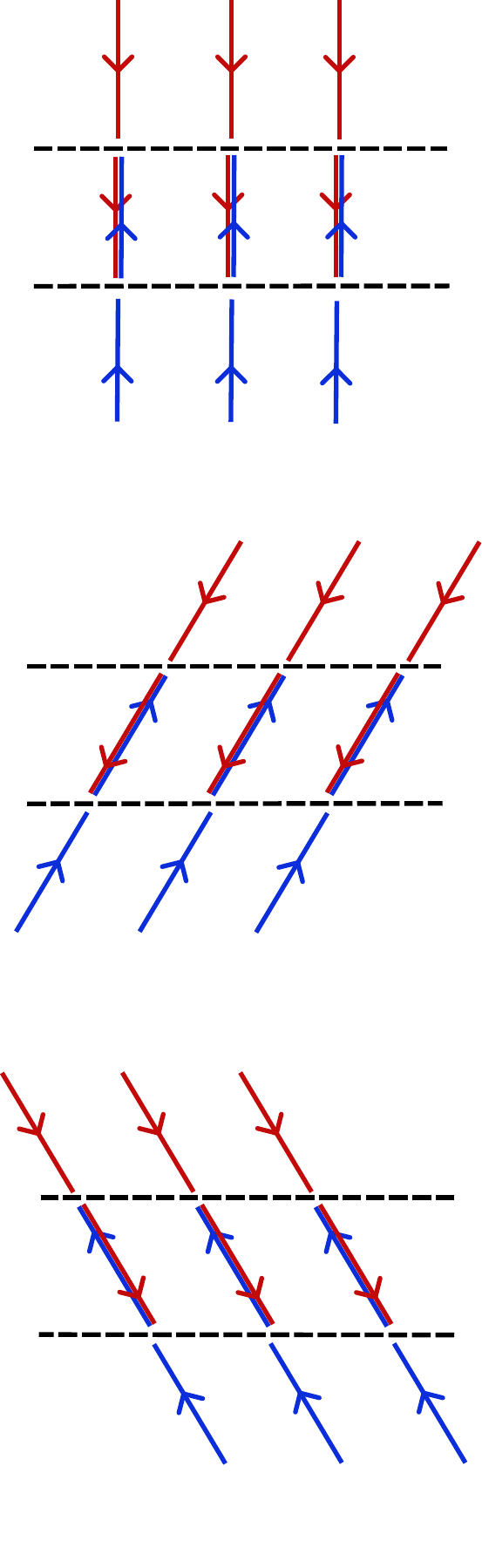}
			\put(10,68){Case $1.27$}
            \put(-3,35){Case $1.58$ with $a_2/b_2=a_1/b_1$}
            \put(-2,1){Case $1.44$ with $a_2/b_2=a_1/b_1$}
		\end{overpic}
	\end{center}
\caption{Possible cases when $\alpha=0$.}\label{figura_6}
\end{figure}

From now on, we consider 
$$
\alpha:=\dfrac{a_2}{b_2}-\dfrac{a_1}{b_1},
$$
with $b_1b_2\neq 0$. Therefore, if $\alpha=0$, the Poincaré first return map presents infinite fixed points located at the line $y=\mu$, i.e., the flow goes from $p_1=\left(x_0^1, \mu\right)$ under the action of $X_1$ until reach the line $y=-\mu$ and comes back under the action of $X_2$ to the same point $p_1$ at the line $y=\mu$, as showed in Figure \ref{figura_6}. We stress that $\alpha=0$ implies that  the vectors $\left(a_1, b_1\right)$ and $\left(a_2, b_2\right)$ are linearly dependent. The cases that represent $\alpha=0$ are the cases 1.44, 1.58 when $a_2/b_2=a_1/b_1$ and the Case 1.27

If $\alpha<0$, taking an arbitrary initial condition $(x_0,y_0)$, the trajectory will converge to a monotone zig-zag behavior restricted to HB region, with the $x$-coordinate going to $- \infty$. The cases that represent $\alpha<0$ are the cases 1.44, 1.58 when $a_2/b_2<a_1/b_1$ and the cases 1.28, 1.59, 1.60. See Figure \ref{figura_1.4}, Figure \ref{figura_1.6} and Figure \ref{figura_1.8}.

If $\alpha>0$, taking an arbitrary initial condition $(x_0,y_0)$, the trajectory will converge to a monotone zig-zag behavior restricted to HB region, with the $x$-coordinate going to $+ \infty$. The cases that represent $\alpha>0$ are the cases 1.44, 1.58 when $a_2/b_2>a_1/b_1$ and the cases 1.26, 1.42, 1.43. See Figure \ref{figura_1.4}, Figure \ref{figura_1.6} and Figure \ref{figura_1.8}.

\section{Conclusion}

In this paper, we presented what is, to the best of our knowledge, the first rigorous mathematical classification of the limit sets arising in a class of hysteretic switching systems frequently encountered in practical applications. Motivated by control protocols in which the switching between different operating regimes occurs through two distinct thresholds, we considered the simplest nontrivial configuration: two planar constant vector fields separated by a hysteresis band bounded by two switching manifolds.

For this class of systems, we completely characterized the asymptotic behavior of trajectories. We showed that every solution exhibits one of three possible qualitative behaviors: unbounded motion, convergence to a monotone zig-zag dynamics confined to the hysteresis band, or periodic behavior. The analysis was performed through a detailed case-by-case study of the vector field configurations and, whenever necessary, through the construction of a Poincaré first return map.

A particularly interesting result is that the existence of periodic orbits and monotone zig-zag dynamics can be determined by a simple parameter involving the slopes of the two vector fields. This provides a clear geometric interpretation of the mechanisms responsible for the long-term behavior of the system.

The present work should be regarded as a first step toward a broader theory of hysteretic switching systems with large hysteresis. Several natural extensions remain open, including the study of affine and nonlinear vector fields, higher-dimensional systems, multiple switching regimes, and more general switching manifolds. We believe that the framework developed here provides a useful foundation for future investigations and for the mathematical understanding of hysteresis-based control protocols arising in medicine, biology, engineering, economics, and other applied sciences.

\section*{Acknowledgments}

Tiago Carvalho is partially supported by S\~{a}o Paulo Research Foundation (FAPESP grants \# 2024/15612-6, \#2021/12395-6, \#2019/10269-3 and \#2022/02819-6) and by Conselho Nacional de Desenvolvimento Cient\'{i}fico e Tecnol\'{o}gico (CNPq Grants 309378/2023-0 and 401974/2025-1).

Bruno de Souza Rangel is supported by S\~{a}o Paulo Research Foundation (FAPESP grant \#2023/18081-9).

\section*{Author Contributions}
All authors contributed equally to this work, including conception, design, analysis, and writing. All authors have read and approved the final manuscript.

\section*{Data availability}
Data sharing is not applicable to this article, as no data sets were generated or analyzed during the current study.

\section*{Declaration of interest}
The authors declare that they have no conflict of interest.

\end{document}